\documentclass [12pt]{article}      
\usepackage{amsmath}
\usepackage{latexsym}
\usepackage{amssymb}
\usepackage{color}
\usepackage{geometry}         
\usepackage{pb-diagram}       
\usepackage{graphicx}

\newtheorem{theorem}{Theorem}
\newtheorem{corollary}[theorem]{Corollary}
\newtheorem{lemma}[theorem]{Lemma}
\newtheorem{definition}[theorem]{Definition}
\newtheorem{prop}[theorem]{Proposition}
\newtheorem{proposition}[theorem]{Proposition}

\newtheorem{ex}[theorem]{Example}
\newtheorem{example}[theorem]{Example}

\newtheorem{remark}[theorem]{Remark}
\newtheorem{claim}[theorem]{Claim}
\newtheorem{exercise}[theorem]{Exercise}

\newcommand{\cat}{^\frown}
\newcommand{\rest}{\ensuremath{\! \upharpoonright\!} }

\newcommand{\sig}{\ensuremath{{\Sigma\!\!\!\!_{_\sim}\, }}}

\newcommand{\omom}{\ensuremath{{\omega^{\omega}}}}

\newcommand{\qed}{{\nopagebreak \hfill $\dashv$  \par\bigskip}}

\newcommand{\pf}{{\par\noindent{$\vdash$\ \ \ }}}

\newcommand{\poQ}{\ensuremath{\mathbb Q}}

\newcommand{\la}{\langle}
\newcommand{\ra}{\rangle}

\newcommand{\mca}{\ensuremath{\mathcal A}}

	  \newcommand{\trees}{\mathcal{T}\!\! rees}

\newcommand{\abs}[1]{\left|#1\right|}

\newcommand{\chisub}[1]{{\chi_{\lower2.5pt\hbox{$\scriptstyle #1$}}}}

\newcommand{\beq}{\begin{equation}}
\newcommand{\eeq}{\end{equation}}

\newcommand{\beqa}{\begin{eqnarray}}
\newcommand{\eeqa}{\end{eqnarray}}

\newcommand{\beqaN}{\begin{eqnarray*}}
\newcommand{\eeqaN}{\end{eqnarray*}}

\renewcommand{\phi}{\varphi}

\newcommand{\bsigma}{\mbox{$\boldmath \Sigma{\boldmath}$}}
\newcommand{\bpi}{\mbox{$\boldmath \Pi{\boldmath}$}}
\newcommand{\bdelta}{\mbox{$\boldmath \Delta{\boldmath}$}}
\renewcommand{\subset}{\subseteq}

\newcommand{\om}{\omega}

\title{Naive Descriptive Set Theory
}
\author{M Foreman}

\begin{document}

\maketitle

%
%






\section{Introduction}
These notes are meant as an introduction to descriptive set theory. The goal is  to provide a propaedeutic tour through the basics of the subject that does not assume any background in Set Theory (well, almost no background). The material is based on some notes taken by Bjork at a course given by Foreman at UC Irvine, and extensively rewritten by Foreman with editing help from the students in a graduate course given in 2010 at the Hebrew University of Jerusalem.
 
Originally intended as an appendix for the paper \emph{The conjugacy problem in ergodic theory}, by Foreman, Rudolph and Weiss, the notes contain adequate background for the descriptive set theoretic portion of that paper and covers some additional topics, such as \emph{norms}, useful for other applications in analysis.

Most of the material is taken from the book \emph{Classical Descriptive Set Theory}, by A. Kechris \cite{Kechbook}, as well as the lecture notes from D. Marker (\cite{marker}) that are on the web.
The notes here are much less complete than either of those sources.

\subsection{Set Theoretic Background and notation}
In this section we briefly describe the Set Theoretic background we will use. It can be safely skipped and referred back to when necessary.

We will avoid the controversy about whether $0$ is a natural number by using the notation $\omega$ for the set $\{0, 1, 2, 3, \dots\}$. People uncomfortable with this notation can systematically substitute $\mathbb N$, with the convention that $\mathbb N$ contains $0$.

We use the von Neumann definition of natural numbers and ordinals. In particular
we identify $0$ with the empty set. Recursively, we identify the number $n$ with the set $\{0, 1, 2, \dots n-1\}$. For example,  the number $4$ is
	the set $\{0,1,2,3\}$.  We note that the number $n$ has exactly $n$ elements.
	
Well-orderings in general, and ordinals in particular play an important role in Descriptive Set Theory. Recall that a well-ordering is a linear ordering $(I,<_I)$ with the property that if $A\subseteq I$ is a non-empty set then $A$ contains an $<_I$ minimal element. Any standard set theory book (such as Levy \cite{levy}) contains the pertinent information.

We summarize the relevant facts:
\begin{itemize} 
\item Each ordinal $\alpha$ is the set of smaller ordinals. The first few ordinals are 
\[0, 1, 2, 3, \dots\] The first infinite ordinal is $\omega$.

\item If $\alpha$ is an ordinal, the least ordinal greater than $\alpha$ is $\alpha\cup \{\alpha\}$ which is denoted $\alpha+1$. Ordinals of the form $\alpha+1$ are \emph{successor} ordinals, the other ordinals are \emph{limit} ordinals.
\item $\alpha$ is an ordinal iff 
\begin{enumerate}
\item $\in$ linearly orders $\alpha$ (i.e. $(\alpha, \in)$ is a linear ordering) and
\item If $x\in \alpha$ then $\in$ linearly orders $x$.

\end{enumerate}	
\item We will write $OR$ for the class of ordinals.
\end{itemize}

By far the most important property of ordinals is that they give canonical examples of well-orderings. This is summarized in the following proposition:
\begin{proposition}
Suppose that $(I, <_I)$ is a well-ordering. Then there is a unique ordinal $\alpha$ and a unique bijection $f:I\to \alpha$ such that 
\[i<_Ij\mbox{ iff } f(i)\in f(j).\]
In other words there is a unique ordinal $\alpha$ such that $(\alpha, \in)$ is isomorphic with $(I,<_I)$ and the isomorphism is unique.

\end{proposition}

In the same constellation of results we have the following definition and proposition.
\begin{definition}\label{well founded relations}
Let $X$ be a set and $R\subseteq X \times X$ be a relation. Then $R$ is \emph{well-founded} iff every 
non-empty subset  $A\subseteq X$ has an $R$-minimal element; i.e. there is an $a\in A$ such that for all $b\in A$, $(b,a)\notin R$. 
\end{definition}

It is worth noting that the following proposition does not use the \emph{Axiom of Choice}.

\begin{prop}
Let $R\subset X\times X$ be a relation. Then $R$ is well-founded iff there is a function $f:X\to OR$ such that for all $a, b\in X$ if $(a,b)\in R$ then $f(a)\in f(b)$.
\end{prop}

We will mostly be interested in countable ordinals. The set of countable ordinals is itself an ordinal. This ordinal is the least uncountable ordinal and will be denoted $\omega_1$.

Other set theoretic notation we will use includes:

\begin{enumerate}
\item We will write $A^B$ for the collection of functions from $B$ to $A$.
\item  \emph{Warning:} Many (most?) texts write $^B\!A$ for the functions from $B$ to $A$.
\item Examples of this include $2^\omega$, the collection of functions from $\omega$ to $\{0,1\}$ and $\omega^\omega$, the collection of functions from $\omega$ to $\omega$. 
\item We will consider \emph{finite sequences} of elements of $B$ to be functions $s:n\to B$ where $n\in \omega$. 
\item The \emph{collection} of finite sequences of elements of $B$ will be denoted $B^{<\omega}$. The collection of sequences of elements of $B$ of length $n$ is denoted $B^n$.
\end{enumerate}

\subsection{Polish Topologies}

\begin{definition}
A topological space $(X,\tau)$ is called \emph{Polish} if $\tau$ is separable and if there is a complete 
metric $d$ on $X$ which generates the topology $\tau$.
\end{definition}

This definition is intended to be sufficiently concrete to be able to prove non-trivial theorems, while being abstract enough to encompass, not only the obvious examples such as $\mathbb R$, $\mathbb C$, and $\mathbb{T}$ (etc.) but also more exotic spaces such as spaces of compact sets.

We explicitly \emph{allow} our spaces to have isolated points. For example $\omega$ with the discrete topology is a Polish space.

 A Polish space without isolated points will be called \emph{perfect}. While the empty set technically fulfills this criterion, we will ignore that problem in practice.
 
We note that being Polish is a property of the topology, not of the metric. So the open unit interval $(0,1)$ is a Polish space (as it is homeomorphic to $\mathbb R$) even though the usual metric on $(0,1)$ is not complete. Indeed Polish spaces can admit many complete separable metrics. The next proposition shows that a Polish space always admits a complete metric bounded by $1$.

\begin{proposition}
Suppose $(X,\tau)$ is a Polish space. Let $d$ be any complete metric on $X$ generating $\tau$. Let
\[d^\prime(x,y)=\frac{d(x,y)}{1+d(x,y)}.\]
Then $d^\prime$ is also a complete metric on $X$ which generates $\tau$. \end{proposition}
The metric $d'$ is sometimes called the ``nearsighted" metric associated with $\tau$.

We will use the standard fact that open sets in a separable metric space are countable unions of closed sets (i.e. $\mathcal F_\sigma$) and, by taking complements, closed sets are countable intersections of open sets (i.e. $\mathcal G_\delta$).

\subsection{Product topologies}

Let $(X_i,\tau_i)$ be a topological space for each $i\in I$. The product topology on $X=\Pi_{i}X_i$ has basic open sets those
 $O=\Pi_{i}O_i$, 
where $O_i=X_i$ for all but finitely many indices $i$, and otherwise $O_i$ is open in $X_i$.
The following is an easy exercise.
\begin{theorem}
Let $X_i$ be a Polish space for all $i\in\omega$. Then $X=\Pi_{i\in\omega}X_i$ is Polish.
\end{theorem}

We can explicitly describe the complete separable metric on $\prod X_i$. Let $d_i$ be a complete metric generating the Polish topology on $X_i$ and bounded by $1$. Define the following metric on $X$:
	\[d(x,y)=\sum_{i\in\omega}\frac{d_i(x_i,y_i)}{2^{i+1}}\]
	Then $d$ is a complete separable metric generating the product topology on $X$. 	

Notice that the product topology is the ``topology of pointwise convergence": given sequences $\vec{x}^n=\la x_i^n:i\in \omega\ra$, then $\vec{x}^n$ converges to $\vec{y}$ iff for each $i, x_i^n$ converges in $X_i$ to $y_i$.

If all of the $X_i$ are the same space $X$ then we can identify $\prod_{i\in I}X_i$ with $X^I$. It is particularly common to identify $\omega^\omega$ with $\prod_{i\in\omega}\omega$ and $2^\omega$ with $\prod_{i\in \omega}2$.

\subsection{The Cantor Space and the Baire Space}

There are two spaces that (along with the discrete space $\omega$) will play a particularly important role for us. This is because of their \emph{universality properties} and because they are very easy to use for \emph{coding} and \emph{diagonal arguments}.

The first such space is the Baire Space. We note that  $\omega^{\omega}$ can be canonically identified with $\Pi_{n\in\omega}\omega$. If we equip $\omega$ with the discrete topology, then the Baire Space is $\omega^\omega$ with the product topology. It is  denoted by $\mathcal{N}$ in some texts. 

The Baire Space is a very familiar space in another guise: ``continued fraction" expansions of irrationals in $(0,1)$ give elements of $\omega^\omega$ and the map sending an irrational to its continued fraction expansion is a homeomorphism between $(0,1)\setminus \mathbb Q$ (with the induced topology from $(0,1)$) and $\omega^\omega$.

\begin{example} Consider the following metric on $\omega^\omega$:
\[ d(x,y)=\left\{ \begin{array}{ll}
			 \frac{1}{n+1} & \mbox{where $n$ is the first coordinate such that $x_n\not=y_n$}\\
						0 & \mbox{ if no such $n$ exists}
			   \end{array}
		\right.\]
Then this metric is complete and generates the product topology on $\omega^{\omega}$.
\end{example}
\noindent This metric will be convenient when we want to discuss the ``radius" of sets in the Baire Space.

We will see later (Theorem \ref{Baire Space surjection}) that there is always a continuous surjection from the Baire Space onto any Polish Space.

The other special space we will consider is the Cantor Space. We will view the Cantor Space as  the set $2^\omega$ consisting of infinite sequences of $0$'s and $1$'s. Viewing $2^\omega$ as the product space $\prod_{i\in\omega}2$ we can equip $2$ with the discrete topology and $2^\omega$  with the product topology. 

Clearly $2^\omega\subseteq \omega^\omega$ and the Cantor Space is identified with a closed compact subset of the Baire Space. 

The Cantor Space is homeomorphic with the usual binary Cantor Set in $(0,1)$.  It is universal in the sense that any (non-empty) Polish space with a perfect subset contains a closed subspace  homeomorphic to The Cantor Space. (See Theorem \ref {Cantor set universality}.)

\begin{remark}The spaces $\omom$ and $2^\omega$  are quite flexible. For example it is easy to verify the following:
\begin{enumerate}
\item Let $b:\omega\to \omega\times \omega$ be a bijection. Then the map $\phi_b:\omega^{\omega \times \omega}\to \omom$ defined by $\phi_b(f)(n)=f(b(n))$ is a homeomorphism that restricts to a homeomorphism of $2^{\om\times \om}$ with $2^\om$.
\item The map $\la f_n:n\in\omega\ra\mapsto g$ defined by $g(n,m)=f_n(m)$ is a homeomorphism between ${(\omom)}^\om=_{def}\prod_{n\in\omega}\omom$ and $\omega^{\omega\times \om}$ that restricts to a homeomorphism of $(2^\om)^\om$ and $2^\om$.

\item  Define a map $\phi:\omom\times \omom\to \omom$ by 
\[\phi(f,g)(k)=\left\{\begin{array}{ll}
			f({k\over 2}) & \mbox{ if  $k$ is even}\\
			g({n-1\over 2}) & \mbox{if $k$ is odd.}
		\end{array}
		\right.\]

Then $\phi$ is a homeomorphism that restricts to a homeomorphism between $2^\om\times 2^\om$ and $2^\om$.
\item The map $(k,f)\mapsto g$ where 
\[g(n)=\left\{\begin{array}{ll}
			k & \mbox{ if }n=0\\
			f(n-1) & \mbox{otherwise}
		\end{array}
		\right.\]
is a homeomorphism between 	$\omega\times \omom$ and $\omom$.
\end{enumerate}
\end{remark}
The point of this remark is that the product topology is given by ``finite information" on a discrete set.
\subsection{Finite Sequences and canonical bases} 
One of the reasons that the Baire Space and the Cantor Space are particularly useful is that their topologies can be understood explicitly as combinatorial objects.

Recall that for any set $A$, we write $A^{<\omega}$ for the set of finite sequences of elements of $A$ viewed as functions $s:\{0,1,2,\dots,k -1\}\rightarrow A$ for some $k$. If we give $A$ the discrete topology then we can identify a base for the product topology on $A^\omega$ as follows.

If $s\in A^{<\omega}$ has length $k$, then we write 
$s=(s_0,s_1,\dots,s_{k-1})$ and $lh(s)=k$. If $a\in A$, we will write $s\cat a$ for the sequence
$\la s_0, \dots s_{k-1}, a\ra$.
If $f\in A^\omega$ we will write $f{\rest} k$ for the sequence $\la f(0), f(1), \dots f(k-1)\ra.$ For $s\in A^{<\omega}$ of length $k$, let $[s]=\{f\in A^\omega\ |\ 
f{\rest} k=s\}$.

\begin{prop}The sets of the form $[s]$ form a clopen basis for the product topology on $A^\omega$.
\end{prop}
\noindent In particular $A^\omega$ is zero-dimensional (totally disconnected). 

For these notes we will be considering the case where $A$ is countable or finite; however it is one of the deep insights of Descriptive Set Theory, in the context of Determinacy, that for uncountable $A$ closed sets in $A^\omega$ code complicated subsets of the Baire Space.

If $A$ is countable or finite, then $A^{<\omega}$ is countable. Hence the basis we have given is countable.  

We will frequently use a combinatorial characterization of convergence in $A^\omega$. From the description of the basis, we see immediately that:
\begin{prop}\label{criterion for convergence}
Let $A$ be a set with the discrete topology. Then $\la x_n:n\in\omega\ra$ converges to $y$ in $A^\omega$ iff for each $k\in \omega$ and all large enough $n, x_n\rest k=y\rest k$.

\end{prop}
\subsection{Infinite trees}
We will analyze closed subsets of $A^\omega$ in terms of \emph{branches through trees}.

\begin{definition}
	Let $A$ be a set. We say that $T\subseteq A^{<\omega}$ is a \emph{tree} if and only if for every $s\in T$, if $s$ has length $k$ and if
	$j\le k$, then $(s_0,\dots, s_{j-1})\in T$. A node $s\in T$ is \emph{terminal} iff it has no proper extension in $T$. A tree $T$ is \emph{pruned} iff it has no terminal nodes.
\end{definition}

\begin{definition}
	Let $T\subseteq A^{<\omega}$ be a tree. A function $f:\omega\rightarrow A$ is an \emph{infinite path} through $T$ iff
	for all $k\in\omega$, $(f(0),f(1),\dots,f(k-1))\in T$. We write $[T]=\{f\ |\ f\mbox{ is an infinite path through T}\}$.

\end{definition}

The next proposition characterizes closed sets as the collections of paths through pruned trees.

\begin{proposition}\label{tree representation of closed sets}
Let $C\subseteq A^\omega$. Then the following are equivalent:
\begin{enumerate}
\item $C$ is closed.
\item There is a tree $T$ such that $C=[T]$.
\item There is a pruned tree $T$ such that $C=[T]$.
\end{enumerate}
\end{proposition}
We note that the pruned tree corresponding to $C$ is canonical.

If $T\subset A^{<\omega}$ is a tree and $s\in T$, then $s$ is a \emph{splitting node} of $T$ iff there are $a\ne b$ in $A$ such that both $s\cat a, s\cat b\in T$.
\begin{proposition}
Let $C$ be a closed set and $T$ a pruned tree with $C=[T]$. Then $C$ is perfect iff for all $t\in T$ there is an $s\in T$ extending $t$ that is a splitting node.
\end{proposition}

\subsection{``Lightface" theory}
There is a recursive enumeration $\la s_n:n\in  \omega\ra$ of finite sequences of natural numbers.  In view of Proposition \ref{tree representation of closed sets}, we can identify closed subsets of the Baire Space with sets of finite sequences of natural numbers, and thus with sets of natural numbers.

By this device we can make sense of such ideas as ``recursive" closed (or open) sets in the Baire Space. By using a recursive pairing function we can make sense of the notions of ``recursive unions" and other natural operations.  Moreover, as we will see, we can identify continuous functions with their neighborhood diagrams and hence we get a well defined notion of a ``recursive" function from $\omega^\omega$ to $\omega^\omega$.  

These ideas turn out to be very useful and important in many contexts, but will not play a significant role in these notes. Hence the title \emph{naive} descriptive set theory.

\subsection{Universality properties of Polish Spaces}
We begin with a discussion that, in some sense, characterizes Polish spaces.
\begin{definition}
The infinite dimensional Hilbert cube is $\mathbb{H}=[0,1]^\omega$ in the product topology, where $[0,1]$ is equipped with the usual
 topology coming from the  real numbers. 
\end{definition}

Viewing elements of $\mathbb H$ as functions $f:\omega\to [0,1]$ standard complete separable metric on the Hilbert Cube is given by 
\[d_H(f,g)=\sum_{i\in\omega}\frac{|f(i)-g(i)|}{2^{i+1}}.\]

\begin{theorem}
	Every Polish space is homeomorphic to a subspace of the Hilbert cube.
	
\end{theorem}
\pf
	Let $X$ be Polish, and $D=\{x_i\ |\ i\in\omega\}$ be any dense subset. Let $d_X$ be a complete metric on $X$ compatible with the 
	Polish topology of $X$ such that $\abs{d_X}\leq1$. Define the following function $f:X\rightarrow\mathbb{H}$ coordinatewise by:
	\[f(x)(i)=d_X (x,x_i)\]
	\begin{claim}\label{claim 1}
		$f$ is uniformly continuous.
	\end{claim}
	\pf[Proof of Claim \ref{claim 1}]
		Let $\varepsilon>0$ and $\delta<\frac{\varepsilon}{2}$.
		Suppose $d_X(x,y)<\delta$. For each $i$, $d_X(x,x_i)\leq d_X(x,y)+d_X(y,x_i)$ and $d_X(y,x_i)\leq d_X(y,x)+d_X(x,x_i)$. Hence
		\[\abs{d_X(x,x_i)-d_X(y,x_i)}\leq d_X(x,y)<\frac{\varepsilon}{2}\]
		It follows that 
		\begin{eqnarray*}
			d_\mathbb{H}(f(x),f(y)) & = & \sum_{i\in\omega}\frac{\abs{d_X(x,x_i)-d_X(y,x_i)}}{2^{i+1}}\\
								& \leq & \sum_{i\in\omega}\frac{\varepsilon}{2^{i+2}}\\
								& < & \varepsilon
		\end{eqnarray*}
	\qed
	\begin{claim}\label{Claim 2}
		$f$ is injective.
	\end{claim}
	\pf[Proof of Claim \ref{Claim 2}]
		Suppose $x\not=y$. By the density of $D$, choose $x_i$ such that $d_X(x,x_i)<\frac{d_X(x,y)}{2}$. Then $d_X(x,x_i)\not=d_X(y,x_i)$ so that
		$f(x)(i)\not=f(y)(i)$.
		\qed
	\begin{claim}\label{Claim 3}
		$f^{-1}$ is  continuous on the range of $f$.
	\end{claim}
	\pf[Proof of Claim \ref{Claim 3}]
		Let $\varepsilon>0$, $x\in X$. Choose $n\in\omega$ such that $d_X(x,x_n)<\frac{\varepsilon}{3}$. Let $\delta=\frac{\varepsilon}{3\cdot2^{n+1}}$.
		Assume, for a proof by
		contraposition, that $d_X(x,y)\geq\varepsilon$ for some $y\in X$. Then the triangle inequality implies that 
		$d_X(y,x_n)\geq\frac{2\varepsilon}{3}$ so that
		
		\begin{eqnarray*}
		d_\mathbb{H}(f(y),f(x)) & \geq & {\abs{d_X(x,x_n)-d_X(y,x_n)}\over 2^{n+1}} \\
					& \ge& \frac{\frac{\varepsilon}{3}}{2^{n+1}}\\
					& = & \delta\end{eqnarray*}

	\qed
	Thus $X$ is homeomorphic to $f(X)\subseteq\mathbb{H}$.
\qed
We will see shortly (Theorem \ref{G delta}) that a subset of a Polish space is Polish with the induced topology iff it is a $\mathcal G_\delta$ set. Once we have that result we can draw the following corollary:
\begin{corollary}
Let $(X,\tau)$ be a topological space. Then $X$ is Polish iff it is homeomorphic to a $\mathcal G_\delta$ subset of the Hilbert cube.
\end{corollary}

We now discuss the universality properties of $\omega^\omega$. We will use these to reduce questions about arbitrary Polish Spaces to questions about $\omom$.
\begin{lemma}\label{F-sigma}
	Let $X$ be a Polish space, $d$ a complete metric, $D\subseteq X$ be an  $\mathcal F_\sigma$ set and $\varepsilon>0$. Then there is a pairwise disjoint sequence of $\mathcal F_\sigma$ sets $\la D_i\ |\ i\in\omega\ra$ such that 
	\[D=\bigcup_{i\in\omega} D_i=\bigcup_{i\in\omega}\overline{D_i}\]
	and $diam(D_i)<\epsilon$ for each $i\in\omega$.

\end{lemma}
\pf Since each closed set can be written as a countable union of closed sets of radius less than $\epsilon$ we can assume that $D=\bigcup_{j\in\omega} D'_j$ where each $D'_j$ has 
diameter less than $\epsilon$. Let $D_i= D'_i\setminus \bigcup_{j<i}D_j=D_i'\setminus \bigcup_{j<i}D_j'$.
 Then each $D_i$ is the intersection of an open set with a closed set and hence is an $\mathcal F_\sigma$ set. 
 \qed

We will frequently be building ``Schemes" in our constructions. These are combinatorial objects usually consisting of trees labelled with sets that have certain properties. Here is our first example.

 \begin{prop}\label{F-sigma scheme} Let $X$ be a Polish space with complete separable metric $d$. There is a sequence of $\mathcal F_\sigma$ sets ${\la}D_{s}\ |\ {s}\in\omega^{<\omega}{\ra}$ with the following properties:
	\begin{enumerate}
	\item $D_\emptyset=X$
	
	\item $diam(D_{s})<\frac{1}{lh({s})}$
	\item $\overline{D_\tau}\subseteq D_{s}$ whenever $\tau$ extends ${s}$
	\item $D_{s}=\bigcup_{i\in\omega}D_{{s}\verb2^2 i}$
	\item If $i\ne j$, then $D_{s\cat i}\cap D_{s\cat j}=\emptyset$
	\end{enumerate}
\end{prop}
 \pf Construct $D_{s}$ by induction on the length of ${s}$. At successor stages, use Lemma \ref{F-sigma}.
\qed

\begin{theorem}\label{cont bij}
Suppose that $X$ is a Polish Space. Then there is a closed set $C\subset \omega^\omega$ and 
a continuous bijection $f:C\to X$.
\end{theorem}

\pf Fix a complete separable metric $d$ on $X$. Let $\la D_i:i\in \omega\ra$ be a scheme as in Lemma \ref{F-sigma scheme} .  Let
\[C=\{a\in \omom:\bigcap_{n\in\omega}D_{a\rest n}\ne\emptyset\}.\]
We claim that $C$ is closed. Let $\la a_n\ra\to b$ with $a_n\in C$.  We need to see that $b\in C$, i.e. that $\bigcap_{n\in\omega}D_{b\rest n}\ne \emptyset$.  By passing to a subsequence we can assume that for all $n$ and all $m\ge n, a_m\rest n= b\rest n$.  

Choose an $x_n\in D_{a_n\rest n}$. Since the $D_{a_n\rest n}$ are nested and have diameters going to zero, the $x_n$ form a Cauchy sequence. Let $y=\lim_n x_n$. Then for all $n$:
\[y\in \overline{D_{a_n\rest n}}=\overline{D_{b\rest n}}.\]
Again by the nesting properties of the scheme, we see that
\[ \bigcap_n D_{b\rest n}=\bigcap_n \overline{D_{b\rest n}}.\]
Hence $y\in \bigcap_n D_{b\rest n}$ and so $b\in C$. 

We now note that since the diameter of $D_{a\rest n}< {1\over n}$, if $a\in C$ there is a unique element $x_a$ of $\bigcap_n D_{a\rest n}$. Let $f:C\to X$ be defined by $a\mapsto x_a$.

We claim that $f$ is a continuous bijection. Let $a\ne a'$ be elements of $C$. Choose an $n$ such that $a\rest n\ne a'\rest n$. Then $D_{a\rest n}\cap D_{a'\rest n}=\emptyset$. Hence $f(a)\ne f(a')$. To see $f$ is onto we do some Boolean algebra. We will frequently use this type of calculation in the sequel so we do it explicitly once here.
\begin{claim}\label{BA calc}
\[\bigcap_n\bigcup_{s\in \omega^n}D_s=\bigcup_{a\in \omom}\bigcap_n D_{a\rest n}.\]

\end{claim}

\pf(Claim \ref{BA calc}) Let $x\in \bigcap_n\bigcup_{s\in \omega^n}D_s$. Using the fact that the $D$'s are nested we can inductively construct a sequence of numbers $a_0, a_1, \dots a_i \dots$ such that for all $n$, $x\in D_{\la a_0, \dots a_{n-1}\ra}$. 
Define $a\in \omom$ by setting $a(i)=a_i$. Then \[x\in \bigcap_{n\in\omega}D_{a\rest n}.\]

On the other hand if $x\in \bigcup_{a\in \omom}\bigcap_n D_{a\rest n}$, we can let $s=a\rest n$ to see that for all $n, x\in \bigcup_{s\in \omega^n}D_s$.\qed

We now see that $f$ is onto: since $D_s=\bigcup_{i\in \omega} D_{s\cat i}$, we can inductively show that for all $n, \bigcup_{s\in \omega^n}D_s=X$.

To see the continuity of $f$, let $a\in\omom$ and $\epsilon>0$. Choose $n$ so large that $1/n<\epsilon$. Since the diameter of $D_{a\rest n}<1/n<\epsilon$, we see that for all $b\in[a\rest n]$ we have $d(a,b)<\epsilon$.\qed

We remark that we cannot prove, in general, that $f$ is a homeomorphism, since closed subspaces of $\omom$ are completely disconnected.

\begin{proposition}\label{retract}
Let $C\subseteq \omom$ be closed. Then there is a continuous map $g:\omom\to C$ such that $g\rest C$ is the identity (i.e. $g$ is a \emph{retract} of $\omom$ to $C$).
\end{proposition}
\pf By Proposition \ref{tree representation of closed sets}  we can find a pruned tree $T$ such that $C=[T]$. Since $T$ is pruned, we can recursively define a map $G:\omega^{<\omega}\to T$ with the following properties:
\begin{enumerate}
\item $lh(s)=lh(G(s))$,
\item if $s\in T$, then $G(s)=s$,
\item for all $i$, there is a $j$ such that $G(s\cat i)=G(s)\cat j$.
\end{enumerate}
By the last clause if $a\in \omom$ then $\la G(a\rest n):n\in\omega\ra$ is a coherent collection of elements of $T$ whose union is an infinite branch through $T$. Hence if we define $g(a)=\bigcup_n G(a\rest n)$, $g$ is a function from $\omom$ to $[T]=C$ which is the identity on $C$. 

We must check that $g$ is continuous. Let $c\in C$. Then a basic open neighborhood of $c$ in $C$ is of the form $[c\rest n]\cap C$ where $[c\rest n]$ is the usual neighborhood in $\omom$. If $g(a)=c$ and $b\in [a\rest n]$, we know that $g(b)\in [c\rest n]\cap C$.\qed

Combining Theorem \ref{cont bij} and Proposition \ref{retract}, we easily see:
\begin{theorem}\label{Baire Space surjection}
	If $X$ is a Polish space, then there is a continuous surjection $\phi:\omega^\omega\rightarrow X$.
\end{theorem}	

We now turn to the Cantor Space. Using very similar ideas to Theorem \ref {Baire Space surjection} we  prove:
\begin{theorem}\label{cantor surjection}
Let $X$ be a compact Polish space. Then there is a  continuous surjection 
\[\phi:2^\omega\to X.\]
\end{theorem}

\pf 
   Inductively build a scheme of closed sets          $\la C_s:s\in\om^{<\om}\ra$ such that 
 \begin{enumerate}
 \item for all $s\in\om^{<\om}, C_s=\bigcup_{n\in\om}C_{s\cat n}$,
\item  for all $s\in\om^{<\om},n\in \omega$, the diameter of $C_{s\cat n}$ is less than $1/ lh(s)$,
\item for all $s\in \om^{<\om}$, the set of $n$ such that $C_{s\cat n}\ne \emptyset$ is finite.
\end{enumerate}
To pass from $C_s$ to $C_{s\cat n}$, consider  a finite cover $\la B_i:i< n=n(s)\ra$ of $C_s$ by open sets of diameter less than $1/ (lh(s)+1)$. Then $C_{s\cat i}=C_s\cap \overline{B_i}$ for $i<n(s)$ and $C_{s \cat i}=\emptyset$ otherwise.

 Let $T=\{s\in\om^{<\om}:C_s\ne\emptyset\}$. Then $T$ is a finitely branching tree.
Define a surjection $p:2^{<\om}\to T$ by induction on $lh(s)$ so that for all $s\in 2^{<\om}$\begin{enumerate}
\item $s\subseteq t$ implies $p(s)\subseteq p(t)$,

 and there are distinct  $t_0, \dots t_{k-1}$ strictly extending $s$ such that

\item 
$[s]=\bigcup_{i< k}[t_i]\subseteq 2^\om$ and
\item $\{p(t_i):i<k\}=\{p(s)\cat n:p(s)\cat n\in T\}$.
\end{enumerate}
Then 
 for all $a\in 2^\om, \bigcup_{k\in\om}p(a\rest k)\in[T]$. Moreover 
if $a\in 2^\om$ the sequence $\la C_{p(a\rest k)}:k\in\om\ra$ is a decreasing sequence of non-empty compact sets. 
As a consequence $\bigcap_{k\in\om}C_{p(a\rest k)}\ne \emptyset$.
Now define $\phi:2^\om\to X$ by
\[\phi(a)=\bigcap_{k\in\om}C_{p(a\rest k)}.\]
\qed

An interested reader will expand the proof of Theorem \ref{Cantor set universality} to show the following result:
\begin{proposition}
Suppose that $X$ is a compact, zero-dimensional, non-empty, perfect Polish space. Then $X$ is homeomorphic to the Cantor Space.
\end{proposition}

\noindent We next inject $2^\om$ into each perfect Polish space.

\begin{theorem}\label{Cantor set universality}
	If $X$ is a non-empty perfect Polish space, then there is a homeomorphism from $2^\omega$ to a subspace of $X$.
\end{theorem}
\pf
	 We build a scheme of sets ${\la}U_s\ |\ s\in2^{<\omega}{\ra}$ with the following properties:
	\begin{enumerate}
	\item $U_\emptyset=X$
	\item $U_s$ is open in $X$ and non-empty.
	\item $diam(U_s)<\frac{1}{{lh(s)}}$
	\item $\overline{U}_{s \verb#^# i}\subseteq U_s$
	\item $U_{s\verb2^2 0}\bigcap U_{s\verb2^2 1}=\emptyset$
	\end{enumerate}
	
	To build this scheme recursively, notice that if $U_{s}$ is defined, then it is a non-empty open set containing no isolated points. So there are 
	points $x\not=y$ in $U_{s}$. Then we can choose open sets $U_{s\verb2^2 0}$ and $U_{s\verb2^2 1}$ that do not meet,  separate
	$x$ from $y$ and have very small diameter.
Since $X$ is complete, for any $a\in 2^\omega$, we have that
	\[\bigcap_{i\in\omega}\overline{U_{a{\rest} i}}=\{x_a\}\]
for some $x_a\in X$. So let $\phi:2^\omega\rightarrow X$ be given by $\phi(a)=x_a$.
	
	To see that $\phi$ is continuous, let $\varepsilon>0$. Let $n\in\omega$ be so large that $\frac{1}{n}<\varepsilon$. Fix $a\in 2^\omega$ and let 
	$s=a{\rest} n$. Then for all $b\in [s]$, we have \[d(\phi(a),\phi(b))<\frac{1}{{lh(s)}}<\varepsilon\]
	
	To see that $\phi$ is injective, suppose that $a\not=b$. Then there is some $n\in\omega$ such that $a_n\not=b_n$. By property (5) above,
	we have that $U_{a{\rest} n+1}\bigcap U_{b{\rest} n+1}=\emptyset$. Hence $\phi(a)\not=\phi(b)$.
	
	Since the Cantor Space is compact, a 1-1 continuous map from $2^\omega$ into $X$ is a homeomorphism onto its range.
	\qed

 Showing that a subset of a  topological space contains a perfect subset says more than that it is uncountable or even that it has power at least $2^{\aleph_0}$, it is a refined version of the Continuum Hypothesis.
\begin{corollary}
Any infinite perfect Polish space has cardinality $2^{\aleph_0}$.
\end{corollary}

\begin{remark}
We note that the range of the function $\phi$ constructed in Theorem \ref {Cantor set universality} is compact, since $\phi$ is a homeomorphism. As a consequence the $\phi$ image of a Borel set is Borel.
\end{remark}

\subsection{Subspaces of Polish spaces}

In this section we will find necessary and sufficient conditions for a subspace of a Polish space to be Polish. We will also give a tool
for refining Polish topologies while still remaining Polish. 

Recall the following definition:
\begin{definition}
Let $(X,d)$ be a metric space and $A\subseteq X$. Then for any $x\in X$, \(d(x,A)=_{def}\inf_{a\in A}\{d(x,a)\}\).
\end{definition}
\begin{definition}
Let $(X,\tau)$ be a topological space and $A\subseteq X$. We write $\tau\rest A$ for the induced topology whose open sets are:
\[U\in\tau\rest A\ \Leftrightarrow\ U=A\cap O\ \mbox{for some $O\in\tau$.}\]
\end{definition}

If $(X, \tau)$ is a Polish space with compatible complete separable metric $d$, and $C\subseteq X$ is closed, then $d\rest (C\times C)$ is a complete separable metric on $C$. In particular the induced topology on a closed subset of a Polish space is Polish. For open sets one must work a bit harder:

\begin{lemma}
	Let $(X,\tau)$ be a Polish space, and $O\subseteq X$ be open. Then $(O,\tau\rest O)$ is Polish.
\end{lemma}
\pf
	Let $d$ be any complete nearsighted metric on $X$ generating $\tau$. Define the following metric on $O$:
	\[\hat{d}(x,y)=d(x,y)+\abs{\frac{1}{d(x,X\setminus O)}-\frac{1}{d(y,X\setminus O)}}\]
	It is left as an exercise to verify that $\hat{d}$ is a metric. 
	
	We need to show that $\hat{d}$ is a complete metric on $O$, and that $\hat{d}$ generates $\tau\rest O$. For the latter, notice that 
	for any $x,y$ we have that $\hat{d}(x,y)\geq d(x,y)$. It follows that for any $x\in O$ and any $\varepsilon>0$, 
	$B_{\hat{d}}(x,\varepsilon)\subseteq B_d (x,\varepsilon)$. So assume that $x\in O$ and let $\varepsilon>0$. We need to find some $\delta>0$
	such that $B_d (x,\delta)\subseteq B_{\hat{d}}(x,\varepsilon)$. Since $O$ is open, there is an $r>0$ such that $B_d(x,r)\subset O$. 
	Choose $0<\delta<r$ such that 
	\[		\delta+\frac{\delta}{r(r-\delta)}<\varepsilon\]
	If $y\in O$ is such that $d(x,y)<\delta$, then $d(y,X\setminus O)>r-\delta$. We then have:
	\begin{eqnarray*}
		\hat{d}(x,y) & = & d(x,y)+\abs{\frac{1}{d(x,X\setminus O)}-\frac{1}{d(y,X\setminus O)}}\\
		            & = & d(x,y)+ \abs{\frac{d(y,X\setminus O)-d(x,X\setminus O)}{d(x,X\setminus O)d(y,X\setminus O)}}\\
					& \leq & \delta+\abs{\frac{d(y,X\setminus O)-d(x,X\setminus O)}{r(r-\delta)}}\\
					& \leq & \delta+\frac{\delta}{r(r-\delta)}\\
					& < & \varepsilon
	\end{eqnarray*}
	
	To see that $\hat{d}$ is complete, let ${\la}x_n\ |\ n\in\omega{\ra}$ be a Cauchy sequence with respect to $\hat{d}$. Since $\hat{d}\geq d$, 
	this sequence is Cauchy with respect to $d$. Since $d$ is complete in $X$, there is an $x\in X$ that this sequence converges to $x$ in the 
	$d$ metric. The metrics $d$ and $\hat{d}$ have the same convergent sequences so it suffices to show that $x\in O$. Since the sequence is Cauchy with respect to $\hat{d}$, we have:
	\[\lim_{n,m\rightarrow\infty}\abs{\frac{1}{d(x_n,X\setminus O)}-\frac{1}{d(x_m,X\setminus O)}}=0\]
	Thus \(\lim_{n\rightarrow\infty}\frac{1}{d(x_n,X\setminus O)}\) exists and is finite. Call this limit $r\geq1$. We now have that
	\[0<\frac{1}{r}=\lim_{n\rightarrow\infty}d(x_n,X\setminus O)=d(x,X\setminus O)\]
	so that $x\in O$.

\qed

\begin{theorem}\label{G delta}
	Let $(X,\tau)$ be a Polish space, and $Y\subseteq X$. Then $(Y,\tau\rest Y)$ is Polish if and only if $Y$ is a $\mathcal G_\delta$ set in $X$.
\end{theorem}

\pf Suppose first that $Y$ is a $\mathcal G_\delta$ set in $X$, say $Y=\bigcap_{n\in\omega} U_n$ for some $U_n$ open in $X$. Without loss of generality we can assume that for all $n, U_{n+1}\subseteq U_n$.
	By the previous Lemma, there is a complete  metric $d_n$ on $U_n$ which generates the subspace topology and is bounded by $1$.
	Define the metric $d:Y\times Y\rightarrow [0,+\infty)$ by $d(x,y)=\sum_{n\in\omega} \frac{d_n(x,y)}{2^n}$. 
	Since $\frac{d_n}{2^n}\leq d$, the topology generated by $d$ is finer than the subspace topology. 
	
To see that they
	are the same topology, let $\varepsilon>0$, $x\in Y$. We need to find an open set $O$ in the subspace topology that is included in $B_d(x,\epsilon)$ and contains $x$.  Choose $N$ so large that 
	\[\sum_{n\geq N}\frac{1}{2^n}<\frac{\varepsilon}{6}\]
	Let $O=\bigcap_{n=1}^{N-1} B_{d_n}(x,\frac{1}{2^n}\frac{\varepsilon}{6})$. Then $O\cap Y$ is an open set in the subspace topology on $Y$. Fix $y\in Y\cap O$. Then
	\begin{eqnarray*}
	d(x,y) & = & \sum_{n\in\omega} \frac{d_n(x,y)}{2^n}\\
		   & \leq & \sum_{n=1}^{N-1} \frac{d_n(x,y)}{2^n}+\frac{\varepsilon}{6}\\
		   & \leq & \frac{1}{6}\sum_{n=1}^{N-1}\frac{\varepsilon}{2^n}+\frac{\varepsilon}{6}\\
		   & < & \frac{\varepsilon}{3}
	\end{eqnarray*}
	Hence $O\subseteq B_d(x,\varepsilon)$.

	To see that $d$ is complete, let ${\la}x_k\ |\ k\in\omega{\ra}$ be a Cauchy sequence with respect to $d$. 
	Then for each $n\in\omega$, ${\la}x_k\ |\ k\in\omega{\ra}$ is a Cauchy sequence with respect to $d_n$. Since all
	of these metrics are compatible with $\tau$, these Cauchy sequences must converge to the same point, say $x$.
	Let $\varepsilon>0$ and $N\in \omega$ be such that $\sum_{n\geq N}\frac{1}{2^n}<\frac{\varepsilon}{3}$. Let $k$ be
	so large that if $1\leq n<N$ and $k^\prime\geq k$ then $d_n(x_{k^\prime},x)<\frac{1}{2^n}\frac{\varepsilon}{3}$. Then
	\begin{eqnarray*}
		d(x_k,x) & = & \sum_{n=1}^{N-1}{\frac{d_n(x_k,x)}{2^{n}}}+\sum_{n\geq N}{\frac{d_n(x_k,x)}{2^{n}}}\\
		         & < & \frac{1}{3}\sum_{n=1}^{N-1}\frac{\varepsilon}{2^n}+\sum_{n\geq N}\frac{1}{2^n}\\
				 & < & \varepsilon
	\end{eqnarray*}
	
	Suppose now that $Y$ is a Polish subspace of $X$. Since closed subsets of a Polish space are Polish, we can
	assume that $X$ is the closure of $Y$. Let ${\la}U_i\ |\ i\in\omega{\ra}$ enumerate a basis for the topology of $X$. 
	Let $d$ be a complete bounded metric on $Y$ generating the subspace topology and $d_X$ be a metric on $X$ generating $\tau$. Fix $x\in Y$ and $\varepsilon>0$. 
	Then there is an $i\in\omega$ such that $x\in U_i$, $U_i$ has $d_X$-diameter less than $\epsilon$ and $U_i\cap Y$ has $d$-diameter less than $\epsilon$. Thus
	\[Y\subseteq \left\{x\in X\ |\ \forall\varepsilon>0\ \exists n\in\omega\ x\in U_n\land diam_d(Y\cap U_n)<\varepsilon\land diam_{d_X}(U_n)<\epsilon\right\}\]
	Call this set $A$. Since
	 \[A=\bigcap_{m\in\omega}\bigcup\{U_n:diam_d(Y\cap U_n)<\frac{1}{m}\land diam_{d_X}(U_n)<\frac{1}{m}\},\]
	 $A$ is a $\mathcal G_\delta$ set. It suffices to show that $A\subseteq Y$. 
	
	Let $x\in A$. For each $m$ choose $n(m)\in\omega$ such that $x\in U_{n(m)}$, $diam_d(Y\cap U_{n(m)})
	<\frac{1}{m}$ and $diam_{d_X}(U_{n(m)})<\frac{1}{m}$. Let $x_m\in U_{n(1)}\cap\dots\cap U_{n(m)}\cap Y$. Such an $x_m$ exists as $X$ is the closure of $Y$ and the intersection of the $U_{n(i)}$'s is non-empty.
	Since the diameter of the $U_{n(m)}$'s approach $0$ as $m$ approaches infinity, the sequence ${\la}x_m\ 
	|\ m\in\omega{\ra}$ is Cauchy with respect to $d$. Hence there is a $z\in Y$ to which this sequence converges. 
	In $X$, this sequence converges to $x$, and so $x=z$, and hence $x\in Y$.
\qed

\section{Borel Sets and the Borel Hierarchy}

\subsection{$\sigma$-Algebras}
We start with a long collection of definitions that can be skipped and referred back to.
\begin{definition}
	Let $X$ be a set. A collection $\mathfrak{A}$ of subsets of $X$ is called a $\sigma-$algebra iff
	\begin{enumerate}
		\item $X\in\mathfrak{A}$
		\item If $A\in\mathfrak{A}$, then $X\setminus A\in\mathfrak{A}$.
		\item If ${\la}A_i\ |\ i\in\omega{\ra}$ is a sequence of sets in $\mathfrak{A}$, then both
		\[\bigcup_{i\in\omega} A_i\in\mathfrak{A}\mbox{ and }\bigcap_{i\in\omega} A_i\in\mathfrak{A}\]
	\end{enumerate}
	We write $(X,\mathfrak{A})$ when talking about the structure of the $\sigma-$algebra $\mathfrak{A}$. 
	\end{definition}
	
We note that to verify that a collection of sets is a $\sigma$-algebra it suffices to show that it is closed under complements and \emph{either} countable unions or countable intersections.	
\begin{definition}
	If $\mathcal{A}$ is any collection of subsets of $X$, we write $\sigma(\mathcal{A})$ for the smallest $\sigma-$algebra
	of sets containing $\mathcal{A}$. 
\end{definition}
 If $\mca$ is a collection of subsets of $X$, then the collection of \emph{all} subsets of $X$, called the  \emph{Power Set} of $X$ or $P(X)$, is a $\sigma$-algebra including $\mathcal A$. We can characterize $\sigma(A)$ ``from the outside" by noting that it is the intersection of all $\sigma$-algebras that include $\mathcal A$. We will get more information by ``building" the $\sigma$-algebra by closing under the operations of countable union and complement.

\begin{definition}	If $(X, \tau)$ is a topological space, the \emph{Borel $\sigma-$algebra} of $X$ is defined to be $\sigma(\tau)$,	the smallest $\sigma-$algebra containing all open sets of $X$. We will sometimes denote this as 
$\mathcal B(\tau)$. The members of $\mathcal B(\tau)$ are called Borel sets. 
\end{definition}

\begin{definition}
	Let $\mathfrak A\subset P(X)$ and $\mathfrak B\subseteq P(Y)$ be $\sigma$-algebras, and $f:X\rightarrow Y$. We say that $f$ is $\mathfrak A$-measurable
	iff for every $B\in\mathfrak{B}$, $f^{-1}(B)\in\mathfrak{A}$.
	A  pair $(Y,\mathfrak{B})$ is a \emph{standard Borel space} iff there is Polish space $(X, \tau)$ and a bimeasurable
	bijection $f:(X,\mathcal B(\tau))\rightarrow (Y,\mathfrak{B})$.
	
\end{definition}
If $f$ is a function with domain a topological space, then $f$ is Borel measurable if it is measurable with respect to Borel $\sigma-$algebras.

\subsection{Some easy facts about Borel sets and functions} 
If $f$ is a function then $f^{-1}$ preserves both finitary and infinitary Boolean operations. As a consequence compositions of Borel measurable functions are Borel measurable, and moreover:

\begin{remark}
	Suppose $(X,\tau)$ and $(Y,\sigma)$ are  topological spaces. Then
	\begin{enumerate}
	\item $f:X\rightarrow Y$ is Borel measurable iff for every $A\in\sigma$, $f^{-1}(A)\in\mathcal B(\tau)$.
	\item If $Y$ has a countable base $\la B_i:i\in\omega\ra$ then $f:X\rightarrow Y$ is Borel measurable iff $f^{-1}(B_i)\in\mathcal B(\tau)$ for every
	$i\in\omega$.
	
	\end{enumerate}
	
\end{remark}
We use this to show:

\begin{proposition}\label{graphborel}
	Suppose $Y$ is a second countable Hausdorff topological space and $X$ is a Hausdorff topological space. 
	If $f:X\rightarrow Y$ is Borel measurable, then $graph(f)=\{(x,y)\in X\times Y\ |\ y=f(x)\}$ is a Borel set 
	in $X\times Y$.
\end{proposition}
\pf
	Let ${\la}A_i\ |\ i\in\omega{\ra}$ be an enumeration of a basis for the topology of $Y$. 
	If $A\subseteq Y$ is open and $B\subseteq X$ is Borel, then $B\times A$ is Borel in $X\times Y$.
	Thus, for each $i\in\omega$, the set $f^{-1}(A_i)\times A_i$ is Borel in $X\times Y$.
	Since $X\times (Y\setminus A_i)$ is closed, we have that 
	\[\bigcap_{i\in\omega}\left(\left(X\times (Y\setminus A_i)\right)\cup \left(f^{-1}(A_i)\times A_i \right)\right)\in\mathcal B({X\times Y})\]
	we claim that this is exactly $graph(f)$.
	
	Suppose $f(x)=y$ and $i\in\omega$. If $y\not\in A_i$, then $(x,y)\in X\times (Y\setminus A_i)$. Otherwise, $x\in f^{-1}(A_i)$.
	
	Conversely, suppose that for all $i\in\omega$, $(x,y)\in \left(X\times (Y\setminus A_i)\right)\cup \left(f^{-1}(A_i)\times A_i \right)$.
	If $y^\prime=f(x)$ for some $y^\prime\not=y$, then there is an $i\in\omega$ such that $y\in A_i$ and $y^\prime\not\in A_i$.
	But then $(x,y)\not\in X\times (Y\setminus A_i)$, so that $(x,y)\in f^{-1}(A_i)\times A_i$. Then $y^\prime=f(x)\in A_i$, a contradiction.
\qed

\subsection{The Borel Hierarchy}

We now give an inductive construction of the Borel sets that begins with the open sets. To keep track of the level of complexity 
of the set, we introduce the following  ``logical" notation:

\begin{definition}\label{borel hi}
Let $(X, \tau)$ be a Polish topological space. The levels of the Borel hierarchy are as follows:
\begin{itemize}
	\item $\bsigma_{1}^{0}$ sets are the open subsets of $X$. 
	\item $\bpi_{1}^{0}$ sets are the closed sets.
	\item Suppose the hierarchy has been defined up to (but not including) the ordinal level $\alpha <\omega_1$. Then $A\in$ $\bsigma_{\alpha}^{0}$ if
	and only if	there is a sequence ${\la}B_{i}\ |\ i\in\omega{\ra}$ of subsets of $X$ with $B_{i}\in\bpi^0_{\beta_i}$ such that for each $i\in\omega$, $\beta_i<\alpha$
	and
	\[A=\bigcup_{i\in\omega} B_i\]
	\item $B\in\bpi_{\alpha}^{0}$ if and only if there is a  subset  of $X$, $A\in$ $\bsigma_{\alpha}^{0}$ such that $B=X\setminus A$.
	\item We write $\bdelta_{\alpha}^{0}=\bsigma_{\alpha}^{0}\cap\bpi_{\alpha}^{0}$.
\end{itemize}
\end{definition}
\noindent We note that this is the ``boldface" theory, indicated by the use of boldface capital Greek letters. In many texts this is written with a ``tilde" (as in $\sig$) to emphasize the difference between the boldface and the lightface theory.

\begin{remark} In more classical language, 
	$\bpi_{2}^{0}$ sets are often called $\mathcal G_\delta$ subsets of a space $X$. 
	 The $\bsigma^0_3$ sets are called the $\mathcal G_{\delta \sigma}$ sets. Similarly, $\bsigma_{2}^{0}$ sets are called $\mathcal F_\sigma$ sets, etc.
\end{remark}

We will also abuse the language in a standard way by using $\bsigma^0_\alpha$ (and $\bpi^0_\alpha$ and $\bdelta^0_\alpha$) both as an adjective and to stand for the collection of sets that are $\bsigma^0_\alpha$ (and $\bpi^0_\alpha$ and $\bdelta^0_\alpha$). We will also frequently use these notations to mean the classes of all  subsets of all Polish spaces that have a given complexity.  If we want to point to the collection of $\bsigma^0_\alpha$ subsets of a particular Polish space $X$ we will sometimes write $\bsigma^0_\alpha(X)$.

The following result justifies  Definition \ref{borel hi} :
\begin{lemma}\label{const of borel}
	For all $0<\alpha<\omega_1$
	we have that $\bsigma_{\alpha}^{0}\cup\bpi_{\alpha}^{0}\subseteq\bdelta_{\alpha+1}^{0}$.
The Borel sets of a Polish space $X$ are precisely $\bigcup_{\alpha<\omega_1}\bsigma_{\alpha}^{0}(X)$.	
\end{lemma}

\pf We first note that  $\bsigma^0_1\subseteq \bsigma^0_{2}$. This is the fact that every open set is a countable union of closed sets, i.e. an $\mathcal F_\sigma$ set. It now follows immediately that $\bsigma^0_\alpha\subseteq \bsigma^0_{\alpha'}$ for all $\alpha<\alpha'<\omega_1$.

By passing to complements we see that $\bpi^0_\alpha\subseteq \bpi^0_{\alpha+1}$. Since it is 
 immediate from the definition that $\bpi^0_\alpha\subset \bsigma^0_{\alpha+1}$ for all $\alpha$, we see that for all $\alpha, \bpi^0_\alpha\subseteq \bdelta^0_{\alpha+1}$. 
Since $\bdelta^0_{\alpha+1}$ is closed under complements, we see that $\bsigma^0_\alpha\subseteq \bdelta^0_{\alpha+1}$. 

From these observations it follows that $\bigcup_{\alpha<\omega_1}\bsigma^0_\alpha$ is closed under complements. To see that it is closed under countable unions, suppose that $A_i\in \bsigma^0_{\alpha_i}$ for each $i\in\omega$, where each $\alpha_i$ is a countable ordinal. Then, since a countable union of countable sets is countable (or weaker: $\omega_1$ is regular), there is an ordinal $\alpha^*$ such that for all $i, \alpha_i+1<\alpha^*$. By the first assertion, $A_i\in \bpi^0_{\alpha_i+1}$, and hence $\bigcup_i A_i\in \bsigma^0_{\alpha^*}$.

We have shown that $\bigcup_{\alpha<\omega_1}\bsigma^0_\alpha$ is a $\sigma$-algebra containing the open sets and hence is a $\sigma$-algebra containing all of the Borel sets. 

On the other hand every member of $\bigcup_{\alpha<\omega_1}\bsigma^0_\alpha$ is built from open sets by taking countable unions and complements, and hence is Borel.\qed

A more refined version of Lemma \ref{const of borel} is given by:

\begin{lemma}\label{closure props of levels}
Let $0<\alpha<\omega_1$.
\begin{enumerate}
	\item $\bsigma_{\alpha}^{0}$ is closed under countable unions and finite intersections. 
	\item $\bpi_{\alpha}^{0}$ is closed under countable intersections and finite unions. 
	\item $\bdelta_{\alpha}^{0}$ is closed under finite intersections and unions, as well as taking complements.
	\item Let $f:X\rightarrow Y$ be a continuous function between Polish spaces. If $A\in\bsigma_{\alpha}^{0}(Y)$ (respectively $\bpi$ or 
	$\bdelta$), then $f^{-1}(A)\in\bsigma_{\alpha}^{0}(X)$ (respectively $\bpi$ or $\bdelta$).
	
\end{enumerate}
\end{lemma}
\pf Inspection of the proof of Lemma \ref{const of borel} yield the first items. To see the last, we note again that the inverses of functions preserve Boolean operations. If $f$ is continuous, then the inverse image of an open set is open. The rest of the claim follows by induction.\qed

\begin{corollary}
	Let $X$ and $Y$ be Polish spaces. If $A\in\bsigma_{\alpha}^{0}(X\times Y)$ (respectively $\bpi$ or $\bdelta$), then for all $x\in X$ the set
	$A_{x}=_{def}\{y\in Y\ |\ (x,y)\in A\} \in\bsigma_{\alpha}^0(Y)$ (respectively $\bpi$ or $\bdelta$).
\end{corollary}

\pf Let $f:Y\to X\times Y$ be given by $f(y)=(x,y)$. Then $A_x=f^{-1}(A)$.\qed

We will often compute the complexity of a Borel set by writing its definition as  a logical formula not involving the implication symbol. In such a formula, $\exists n\in\omega$ corresponds to a countable union and $\forall n\in\omega$ corresponds to a countable intersection.

A pitfall of this  method is that the implication symbol in
\[(\exists n\in \omega) (\phi)\implies (\exists n\in\omega) (\psi)\] must be rewritten as 
\[\large(\neg(\exists n\in \omega) (\phi) \lor (\exists n\in\omega) (\psi)\large)\] before its complexity can be recognized.
The complexity is that of the formula
\[ (\forall n\in\omega)(\neg\phi)\lor (\exists n\in\omega)(\psi).\]
A clear understanding of this is not essential at this point; we will illustrate the relevance of the logical notation mostly with examples.

\begin{ex} \begin{enumerate}

\item Let $X$ be a Polish space and $F\subseteq X$ be countable, then $F\in\bsigma^{0}_{2}(X)$.

\item Let $A\subseteq\omega^\omega$ be the set of functions that are eventually injective. Then 
$f\in A$ if and only if 
\[\exists n\forall m_1,m_2>n\left(f(m_1)\not=f(m_2)\lor m_1=m_2\right).\]
This can be written as:
\[\bigcup_{n\in\omega}\bigcap_{{m_1,m_2>n}\atop{m_1\ne m_2}}\{f: f(m_1)\ne f(m_2)\}.\]

Thus $A$ is a $\bsigma_{2}^{0}$ set.
\end{enumerate}
\end{ex}
An important example involves the group of  permutations of the natural numbers denoted as $S_{\infty}$:

\begin{example}\label{S infinity is G delta} We can view $S_{\infty}$ as a subset of
the Baire Space, since $f\in S_{\infty}$ if and only if $f$ is a bijection from $\omega$ to $\omega$. 
Now the set of injective functions in $\omega^\omega$ is closed: $f$ is injective if and only if $\forall m_1,m_2\left(f(m_1)
\not=f(m_2)\lor m_1=m_2\right)$. Thus this is a countable intersection of closed sets of the form 
\[\{f:f(m_1)\ne f(m_2)\}\]

 Furthermore, $f$ is onto if and only if $\forall n\exists m(f(m)=n)$, hence a countable intersection of open sets of the form $\{f:\exists m(f(m)=n)\}$. This implies that the set of onto functions is
a $\bpi_{2}^{0}$ set. $S_{\infty}$ is therefore the intersection of a closed set with a $\bpi_{2}^{0}$ set, which is $\bpi_{2}^{0}$, by Lemmas \ref {const of borel} and \ref {closure props of levels}.
\end{example}

\begin{exercise}
\begin{enumerate}
\item Show that $S_{\infty}$ is a Polish group.
\item Show that $S_\infty$ is not a $\bsigma_{2}^{0}$ subset of $\omom$.

\noindent (Hint: Use the Baire Category Theorem.)
\end{enumerate}
\end{exercise}
\subsection{Examples of coding}
Here are some important examples of coding:
\begin{example}\label{coding}
Let $R\subset \omega\times \omega$. Then we can associate an element of $2^{\omega\times\omega}$ to $R$ by considering its characteristic function $\chi_R:\omega\times \omega\to 2$ defined by  

\[\chi_R(m,n)=\left\{\begin{array}{ll}
			1 & \mbox{ if  }(m,n)\in R\\
			0 & \mbox{otherwise.}
		\end{array}
		\right.\]
If $f\in 2^{\omega\times \omega}$, we will say that $f$ \emph{codes} $R$ if $R=\{(m,n):f(m,n)=1\}$. Then 

\[\mathcal{LO}=_{def}\{f\in 2^{\omega\times \omega}:f\mbox{ codes a linear ordering of a subset of }\omega\}\]
is a closed set.
\end{example}

\begin{example}\label{space of infinite trees}
Let $\la \sigma_n:n\in\omega\ra$ be a 1-1 enumeration of $\om^{<\om}$ such that $\sigma_m\subseteq \sigma_n$ implies $m\le n$. (Shorter sequences come first.) For $T\subseteq \om^{<\om}$, let $f_T\in 2^\om$ be the characteristic function of $\{n:\sigma_n\in T\}$. Then:
\[\{f_T:T\mbox{ is an infinite tree}\}\]
is a $\bpi^0_2$ subset of $2^\om$ 
\end{example}
\begin{exercise}
What is the complexity of \[\{f_T:T \mbox{ is a finitely branching tree}\}?\]
\[\{f_T:T\mbox{ is a tree}\}?\]
\end{exercise}

\noindent We will refer to the  space of trees with the topology induced from $2^\om$ as $\trees$.
\medskip

Using the fact that $2^\om\cong \prod_{n\in\om}(2^{\om^n})^\om$, we can similarly code any structure in a countable first order language. For an $EC_\Delta$ class of structures, the set of codes is a Borel set. To give a simple example:
\begin{example}
Let $G=\la g_n:n\in\om\ra$ be a countable group with $g_0=e_G$. Associate to $G$ an element  $\chi_G\in 2^{\om\times \om \times \om}$ by setting $\chi_G(m,n,p)=1$ iff $g_m\cdot_G g_n=g_p$. Then the collection of $\chi_G$ for $G$ a countable group is a $\bpi^0_2$ set.
\end{example}
Note that we can let $S_\infty$ act on $2^{\om\times\om\times\om}$ by setting $(g\cdot x)(k,l,m)=x(g^{-1}k, g^{-1}l, g^{-1}m)$ for $g\in S_\infty$ and $x\in 2^{\om\times \om\times \om}$. Then
for countable groups $G, H$ we have $G\cong H$ iff there is a $g\in S_\infty$,
\[g\chi_G=\chi_H.\]
More generally this $S_\infty$ action on $2^\om\cong \prod_{n\in\om}(2^{\om^n})^\om$ codes isomorphism for any class of countable structures, a fact that is relevant to model theory.

\begin{example}\label{embedding baire into Cantor}
Let $\phi:\omom\to 2^{\om\times \om}$ be defined by $f\mapsto \chi_f$, where $\chi_f(m,n)=1$ iff $f(m)=n$. Then $\phi$ is continuous and 1-1. Moreover $\phi^{-1}$ is continuous on the range of $\phi$. In fact the range is $\mathcal G_\delta$. This can be seen using Theorem \ref{G delta}, or by the following computation:

If $a\in 2^{\om \times \om}$ then $a$ is in the range of $\phi$ iff
\begin{enumerate}
\item $(\forall m)(\exists n)a(m,n)= 1$
\item $(\forall m)(\forall n)(\forall p)((a(m,n)=1\land a(m,p)=1)\rightarrow n=p)$
\end{enumerate}
We note that this shows that if $X$ is perfect Polish then there is a $\mathcal G_\delta$ subset of $X$ that is homeomorphic to $\omom$.
\end{example}

The next exercise will be useful when we consider analytic and co-analytics sets:

\begin{exercise}\label{borel bijections one}
Let $X$ be a perfect Polish space. Then there is a bijection $f:X\to \omom$ such that 
$f$ and $f^{-1}$ are Borel measurable. (Hint: Build Borel injections in each direction and then use the Cantor-Schroeder-Bernstein theorem.) 
\end{exercise}

\subsection{Universal Sets}
\begin{definition}
Let $X$, $Y$ be any spaces and $\Gamma$ be a point class of $X$, that is, $\Gamma$ is a collection of subsets of $X$. Let $A\subseteq Y\times X$.
Then $A$ is universal for $\Gamma$ if and only if $\Gamma=\{A_y\ |\ y\in Y\}$. 
\end{definition}

\begin{theorem}
	Let $X$ be a second countable metric space. For $1\leq\alpha<\omega_1$ there is a $U_\alpha\in\bsigma_{\alpha}^{0}(2^{\omega}\times X)$ that is universal
	for $\bsigma_{\alpha}^{0}(X)$. 
\end{theorem}

\pf
	We proceed by induction on $\alpha$. 
	Let ${\la}B_i\ |\ i\in\omega{\ra}$ be a basis for the topology of $X$. Let \[U_1=\bigcup_{n\in\omega}\{(f,x)\ |\ f(n)=1\land x\in B_n\}\]
	It is clear from the definition that $U_1$ is open. To see that it is universal, let $O\subseteq X$ be open. Define the following
	function $f:\omega\rightarrow 2$ by $f(n)=1$ if and only if $B_n\subseteq O$. We claim that $O=\left(U_1\right)_{f}$. 
	
	If $x\in\left(U_1\right)_{f}$, then there is an $n\in\omega$ such that $f(n)=1$ and $x\in B_n$. Hence $B_n\subseteq O$, and thus $x\in O$.
	Conversely, since $O$ is open, if $x\in O$  there is an $i\in \omega$ such that $x\in B_i\subset O$. But then $f(i)=1$ and $x\in B_i$,
	so that $x\in \left(U_1\right)_{f}$. 
	
	Now suppose that for each $\beta<\alpha$ we have a $U_\beta\in\bsigma_{\beta}^{0}(2^{\omega}\times X)$ which is universal for $\bsigma_{\beta}^{0}(X)$.
	Let $V_{\beta}=(2^{\omega}\times X) \setminus U_{\beta}$. It is easy to see that $V_\beta\in\bpi_{\beta}^{0}(2^{\omega}\times X)$ and is universal for
	$\bpi_{\beta}^{0}(X)$. If $\alpha$ is a limit ordinal, let ${\la}\beta_i\ |\ i\in\omega{\ra}$ be an increasing sequence of ordinals approaching $\alpha$. 
	If $\alpha$ is a successor ordinal, say  $\alpha=\beta+1$, set $\beta_i=\beta$ for all $i\in\omega$. 
	
	Fix any homeomorphism $F:2^\omega\rightarrow (2^\omega)^\omega$. For $f\in 2^\omega$, denote $F(f)= {\la}(f)_i\ |\ i\in\omega{\ra}$. Define $U_\alpha$ to be:
	\[U_\alpha = \{(f,x)\ |\ x\in\bigcup_{i\in\omega}\left(V_{\beta_i}\right)_{(f)_i}\}.\]
	
We first claim that $U_\alpha$ is $\bsigma^0_\alpha$.
 Since \[\{(f,x)\ |\ x\in\bigcup_{i\in\omega}\left(V_{\beta_i}\right)_{(f)_i}\}	=\bigcup_{i\in\omega}\left\{(f,x)\ 
|\ x\in\left(V_{\beta_i}\right)_{(f)_i}\right\}\] 
	it suffices to show that for a fixed $i$, $B_i=_{def}\{(f, x):x\in (V_{\beta_i})_{(f)_i}\}$ is $\bpi^0_{\beta_i}$. 
	
Let $Y=(2^\omega)^\omega\times X$ and $G:2^\omega\times X\to Y$ be
given by $(f, x)\mapsto (\la (f)_i\ra, x)$. Then for each $i$,  $\{(\vec{f}, x)\in Y:x\in (V_{\beta_i})_{(f)_i}\}\in \bpi^0_{\beta_i}$. Hence $G^{-1}(\{(\vec{f}, x)\in Y:x\in (V_{\beta_i})_{(f)_i}\})=B_i\in \bpi^0_{\beta_i}$.

It remains to show that $U_\alpha$ is universal.
Let $A\in \bsigma_{\alpha}^0(X)$. Then there is a sequence of sets $\la B_i:i\in \omega\ra$ such that for some $\delta_i<\alpha, B_i\in \bpi^0_{\delta_i}$ and $A=\bigcup_i B_i$. Since the sequence $\la \beta_i\ra$ is not bounded in 
$\alpha$ and the $\bpi^0_\xi$'s increase with $\xi$,  we can assume that $\delta_i=\beta_i$. Choose $f_i\in 2^\omega$ such that $B_i=(V_{\beta_i})_{f_i}$. Let $f\in 2^\omega$ be such that $(f)_i=f_i$. 
We see that 
 \begin{eqnarray*}
		x\in(U_{\alpha})_f & \Leftrightarrow & \exists i\ x\in(V_{\beta_i})_{f_i}\\
		                  & \Leftrightarrow & \exists i\ x\in B_i\\
						  & \Leftrightarrow & x\in A\\
	\end{eqnarray*}
Thus $A=(U_\alpha)_f$.
	\qed

 \begin{corollary}
 Let $X$ be a perfect Polish space. Then for all $0<\alpha<\omega_1$ there is a 
 set $U_\alpha\in \bsigma^0_\alpha(X\times X)$ that is universal for $\bsigma^0_\alpha(X)$.
 \end{corollary}
 
 \pf Let $\phi:2^\omega\to X$ be a homeomorphism of $2^\omega$ into $X$. Then the range of $\phi$ is compact, and hence closed. If we let $F:2^\omega\times X\to X\times X$ be defined by setting $f(a,x)=(\phi(a),x)$, then it is easy to check that $F$ takes closed sets to closed sets. Let $V'_1$ be  a $\bpi^0_1(2^\omega\times X)$
set that is universal for $\bpi^0_1$ subsets of $X$. Then $V_1=F[V'_1]$ is $\bpi^0_1$ and is universal for $\bpi^0_1$ subsets of $X$.
If we let $U_1=(X\times X)\setminus V_1$, then $U_1$ is $\bsigma^0_1(X\times X)$ and universal for $\bsigma^0_1$ subsets of $X$.

For $\alpha>1$, life is easier: Let $U'_\alpha\in \bsigma^0_\alpha(2^\omega\times X)$
be universal for $\bsigma^0_\alpha(X)$. Then  $U_\alpha=_{def}F[U'_\alpha]$ is the intersection of a $\bsigma^0_\alpha(X\times X)$ set with the range of $F$. Since the range of $F$ is closed and $\alpha\ge 2$, this is $\bsigma^0_\alpha$ and is universal for $\bsigma^0_\alpha$ subsets of $X$.\qed

\noindent The reader ambitious enough to try to directly construct universal subsets of $X\times X$ for an arbitrary $X$ will begin to appreciate the utility of working directly with $2^\omega$.

The next argument is Cantor's classical diagonal argument cast in this setting.

\begin{corollary} Let $X$ be a (non-empty) perfect Polish space. Then 
	 \[\left\{\bsigma_{\alpha}^{0}(X)\ |\ 1\le\alpha<\omega_1\right\}\]
 is a strictly increasing collection of subsets of $P(X)$.
\end{corollary}

\pf
	Assume towards a contradiction that $\bsigma_{\alpha+1}^{0}(X)\subseteq\bsigma_{\alpha}^{0}(X)$. Then $\bpi^0_{\alpha+1}(X)\subseteq \bpi^0_\alpha(X)$. Consider  a  set $U_{\alpha +1}\in \bsigma^0_{\alpha+1}(X\times X)$ universal for $\bsigma^0_{\alpha+1}(X)$ and let 
	\[D=\left\{x\ |\ x\not\in(U_{\alpha+1})_x\right\}.\] Then $D\in\bpi_{\alpha+1}^{0}$. By our assumption, we have:
	\[\bpi_{\alpha+1}^{0}\subseteq\bpi_{\alpha}^{0}\subseteq\bsigma_{\alpha+1}^{0}\]
so that $D\in\bsigma_{\alpha+1}^{0}$. By the universality of $U_{\alpha+1}$ there is a $y$ such that $D=(U_{\alpha+1})_y$. 
	
Following Cantor:
\[y\in D \iff (y,y)\in U_{\alpha+1}\iff y\notin D.\]
\qed

\begin{exercise}Let $X$ be perfect Polish. 
\begin{enumerate}
\item Deduce that for all $0<\alpha<\omega_1, \bsigma^0_{\alpha+1}(X)\nsubseteq\bsigma^0_\alpha(X)$ directly from the statement that $\bsigma^0_{\alpha+1}(2^\omega)\nsubseteq \bsigma^0_\alpha(2^\omega)$.
\item Show that for all $0<\alpha<\omega_1$ there is no $\bdelta^0_\alpha(X\times X)$ that is universal for $\bdelta^0_\alpha(X)$.
\end{enumerate}
\end{exercise}

\subsection{Changing Topologies}
For the proof of the next theorem, we need the following Lemma:

\begin{lemma}\label{increasing topologies}
Let $(X,\tau)$ be a Polish space, and suppose that ${\la}\tau_n\ |\ n\in\omega{\ra}$ is a sequence of Polish topologies on
 $X$ with $\tau\subseteq\tau_n$ for each $n$. Then the topology $\tau_\infty$ which is generated by $\cup\tau_n$ is Polish.
Moreover, if $\mathcal{B}(\tau_n)=\mathcal{B}(\tau)$ for each $n$, then $\mathcal{B}(\tau_\infty)=\mathcal{B}(\tau)$.
\end{lemma}
	
	\pf
	Write $X_n$ for the space $(X,\tau_n)$ and $X_\infty$ for the product space $\prod_{n\in\omega}X_n$. By our
	earlier results, $X_\infty$ is a Polish space. Let $\phi:X\rightarrow X_\infty$ be defined by
	$\phi(z)=(z,z,z,\dots)$. 
	
	We claim that $\phi[X]$ is closed in $X_\infty$.
To see this, let $y\in X_\infty$ with $y\not\in\phi[X]$. There are then $n,m\in\omega$ with $y(n)\not=y(m)$. Choose open
		sets $U_n,U_m$ such that $y(n)\in U_n$, $y(m)\in U_m$, and $U_n\cap U_m=\emptyset$. Let $V\subseteq X_\infty$
		be the open set $V=X\times \dots\times X\times U_n\times X\times\dots\times X\times U_m\times X\times\dots$, where
		$U_i$ occurs in the $i$th place.
		Then $V$ is an open neighborhood of $y$ that is disjoint from $\phi[X]$.

	Since $\phi[X]$ is closed,  $\phi[X]$ is a Polish space in the induced topology $\nu$ from the product space. The lemma follows if we can show that $\nu$ coincides with the 
	topology given by $\cup_{n}\tau_n$.
	Let $O$ be a basic open set for $\nu$. Then there are $O_1\in\tau_1$, $O_2\in\tau_2$, $\dots$, $O_n\in\tau_n$
	for some $n\in\omega$ such that $O=\{x\in X\ |\ \phi(x)\in O_1\times O_2\times\dots\times O_n\times X\times\dots\}$.
	Thus $O=O_1\cap O_2\cap\dots\cap O_n$ and hence $O$ is in the topology generated by the $\tau_n$'s. 	\qed
\begin{exercise}
Show that in the previous Lemma, if $\mathcal{B}(\tau_n)=\mathcal{B}(\tau)$ for each $n$, then $\mathcal{B}(\tau_\infty)=\mathcal{B}(\tau)$.
\end{exercise}
\begin{theorem}\label{retopologize}
Suppose that $(X,\tau)$ is Polish and $A\subseteq X$ is Borel. Then there is a Polish topology $\tau_A\supseteq
\tau\cup\{A\}$ that has the same Borel sets as $\tau$.
\end{theorem}
	\pf
	We  first consider the case where $A$ is a closed set $F$. Let $S_1=(F,\tau\rest F)$ and $S_2=(X\setminus F, 
	\tau\rest {X\setminus F})$, which are both Polish spaces. Let $d_1$ be a complete compatible  
	metric on $S_1$ and $d_2$ a complete compatible metric on $S_2$, both bounded by $1$. 
	Define \[d_3(x,y)=\left\{\begin{array}{ll}
							d_1(x,y) & \mbox{if both $x,y\in F$}\\
							d_2(x,y) & \mbox{if both $x,y\in X\setminus F$}\\
							2 & \mbox{otherwise}. 
							\end{array}
							\right.\]
	This metric is Polish on $S_1\oplus S_2$ and generates a topology on $X$ that contains both $\tau $ and $F$. Thus for any Polish space $(X, \tau)$ and any closed set $F$ we can extend the topology $\tau$ to contain $F$, be Polish and have the same Borel sets as $\tau$.
	
	Let $G=\{A\subseteq X\ |${ there is a Polish topology }$\tau_A \supseteq\tau$ with the same Borel sets as
	$\tau$ with $A\in \tau_A\}$. We want to show that every Borel set belongs to $G$. Since $G$ contains
	all the open and closed sets, we need only show that $G$ is closed under countable unions and complementation.
	
	Suppose $A\in G$. Then $(X,\tau_A)$ is Polish with $A\in\tau_A$. So $X\setminus A$ is $\tau_A$-closed, and by the first paragraph of the proof we can   extend $\tau_A$ to $\tau_{A,X\setminus A}$. 
	$G$ is therefore closed under taking complements. 
	Suppose ${\la}A_n\ |\ n\in\omega{\ra} $ is a sequence in $G$. Then $\tau_{A_n}$ satisfy the requirements of
	Lemma \ref {increasing topologies}, hence $\cup_{n\in\omega} A_n\in G$ with Polish topology $\tau_\infty$ on $X$.
	\qed
	
	\begin{corollary}\label{borel sets are polish spaces}
		If $(X,\tau)$ is Polish and $Y\subseteq X$ is Borel, then $(Y,\{B\cap Y\ |\ B\in\mathcal{B}(\tau)\})$ is
		a standard Borel space. 
	\end{corollary}

The next example is cautionary:

\begin{example}
Let $X=[0,1]$ and $\tau$ be the usual topology.  Then $A=_{def}\poQ\cap X$ is an $\mathcal F_\sigma$-set.  Let $\sigma$ be the topology generated by $A$ and $\tau$. We claim that $\sigma$ is \emph{not} Polish. 

For if it were, the Baire Category theorem would apply to this topology.  If we write $A=\bigcup_{q\in A}\{q\}$, then, as each $\{q\}$ is a closed set, one of the $\{q\}$ must have non-empty interior. The only way this could happen is if $\{q\}$ is an open set. But an easy argument shows that there are no singleton open sets in $\sigma$.

\end{example}	
	
\section{Analytic Sets}

\begin{definition}
Let $X$ be a Polish space. A set $A\subseteq X$ is analytic if and only if it is the continuous image of a Polish space $Y$. 

\end{definition}

If $Z=X\times Y$ then $\Pi_X:Z\to X$ is the map $(x,y)\mapsto x$ and $\Pi_Y$ is $(x,y)\mapsto y$. Note that if $A\subseteq Z, \Pi_X[A]=\{x:(\exists y)\ (x,y)\in A\}$

We will denote the collection of analytic subsets of $X$ by $\bsigma^1_1(X)$, with the usual abuses of language.
There are many equivalent ways of thinking of analytic sets. The following Proposition gives several of them:
\begin{proposition}\label{laundry list}
	The following are equivalent:
	\begin{enumerate}
	\item $A$ is analytic.
	\item $A$ is the continuous image of $\omega^\omega$.
	\item There is a Polish space $Y$ and a closed set $F\subseteq X\times Y$ such that 
	\[A=\Pi_X(F)\]
	\item There is a closed set $F\subseteq X\times \omega^\omega$ such that 
	\[A=\Pi_X(F)\]
	\item There is a Borel set $B\subseteq X\times \omega^\omega$ such that 
	\[A=\Pi_X(B)\]
	
	\end{enumerate}
	
\end{proposition}
\pf
	The equivalence of the first two conditions is immediate because $\omom$ is Polish and
	every Polish space $Y$ is the continuous image of $\omom$.
		
	For the proof of ($1$) iff ($3$), assume $A$ is analytic and $f:Y\rightarrow X$
	is continuous with $Y$ Polish and $A=rng(f)$.  
	Then the inverse graph of $f \subseteq X\times Y$ is closed. Moreover, 
	\begin{eqnarray*}
		x\in A & \Leftrightarrow & (\exists y)f(y)=x\\
		       & \Leftrightarrow & x\in\Pi_X(graph^{-1}(f))
	\end{eqnarray*}
	
	Conversely, suppose $F\subseteq X\times Y$ is closed. Then $F$ is a Polish space with the induced topology. 
	The function $p: F\rightarrow X$ defined by $p((x,y))=x$ is continuous with $A= rng(p)$ as desired.
	
	The proof that (2) implies (4) is exactly parallel to the proof that (1) implies (3). Moreover it is obvious that (4) implies (5). To see that (5) implies (1) we fix a Borel set $B\subseteq X\times  \omom$ with $A=\Pi_X(B)$.  By Theorem \ref{retopologize} we can extend the topology on $X\times \omom$ to a larger topology making a space $Y$ with underlying set $X\times \omom$ in which $B$ is closed. Since the topology on $Y$ is finer than the original topology on $X\times \omom$ the map $\Pi_X:Y\to X$ is still continuous. Since $B$ is closed, it is Polish with the induced topology. Hence $\Pi_X:B\to X$ is a continuous map from a Polish space to $X$. Finally we note that $A$ is the range of $\Pi_X\rest B$. 
	\qed

From Theorem \ref {retopologize}, it follows immediately that:
\begin{corollary}
If $(X,\tau)$ is Polish and $B\subseteq X$ is Borel then $B$ is analytic. Moreover continuous images of Borel sets are analytic. 
\end{corollary}	
\pf  The second statement is stronger. Let $f:Y\to X$ and $A=f(B)$. By Theorem \ref {retopologize}, we can retopologize $X$ with a $\sigma$ so that $B$ is open and   $\sigma{\rest}B$ is Polish. Since we have added open sets $f\rest B$ is still continuous. Hence $A$ is the image of a Polish space under a continuous function.\qed

\begin{example}
In the notation of Example \ref{space of infinite trees}, show that 
\[\{f_T:T\subseteq \omega^{<\om}\mbox{is a tree with an infinite branch}\}\]
is a $\bsigma^1_1$ subset of $2^\om$. Why does this differ from Example \ref{space of infinite trees}? 
\end{example}

\begin{example}
In the notation of Example \ref{coding}, show that 
\[\{f\in \mathcal{LO}:f\mbox{ is not a well-ordering}\}\]
is a $\bsigma^1_1$-set.
\end{example}

\subsection{Universal analytic sets and a non-Borel analytic set}

\begin{theorem}
	For any Polish space $X$ there is an analytic set $U\subseteq \omega^\omega\times X$ such that, if $A\subseteq X$
	is analytic, then there is an $f\in\omega^\omega$ with $A=U_f=\{x\ |\ (f,x)\in U\}$.
	
\end{theorem}
In other words, there is an analytic set that is universal for the class of analytic sets of a Polish space.

\pf
	Let $V\subseteq \omega^\omega\times (\omega^\omega\times X)$ be closed and universal for closed sets in 
	$(\omega^\omega\times X)$. Let $U=\{(f,x)\ |\ \exists g\in\omega^\omega\ (f,g,x)\in V\}$, which is analytic.
	To show that $U$ is universal, let $A\subseteq X$ be analytic. Then there is a closed 
	$F\subseteq\omega^\omega\times X$ such that $A=\{x\in X\ |\ \exists g\in\omega^\omega\ (g,x)\in F\}$.
	Since $F$ is closed, there is an $f\in\omega^\omega$ such that $F= V_f$. But then
	\begin{eqnarray*}
	U_f & = & \{x\ |\ \exists g\ (f,g,x)\in V\}\\
	    & = &  \{x\ |\ \exists g\ (g,x)\in F\}\\
		& = & A
	\end{eqnarray*}

\qed

\begin{exercise}
Use Example \ref{embedding baire into Cantor} to show that if $X$ is a perfect Polish space, then there is a set $A\subseteq X\times X$  that is universal for analytic subsets of $X$.
\end{exercise}

%

\begin{definition}
Let $\la X_n:n\in\om\ra$ be a sequence of Polish spaces and assume that $\la d_n\ra$ is a sequence of complete metrics witnessing that the $X_n$'s are Polish that are bounded by $1$. Let $Y=\bigcup_{n\in\omega}X_n$ and define a metric on $Y$ by setting:
\[d_Y(x,y)=\left\{\begin{array}{ll}
							d_n(x,y) & \mbox{if both $x,y\in X_n$ for some $n$}\\
							2 & \mbox{otherwise}. 
							\end{array}
							\right.\]
Then the \emph{separated sum} of the $X_n$'s  is the space $(Y, \tau)$ where $\tau$ is the topology of $d_Y$.
\end{definition}

\begin{theorem}
	Let $X$ be a Polish space, and let ${\bsigma^1_1}(X)=\{A\subseteq X\ |\ A\mbox{ is analytic}\}$. Then
 ${\bsigma^1_1}(X)$ is closed under countable unions and intersections. 
	\end{theorem}
\pf
Let ${\la}A_n\ |\ n\in\omega{\ra}$ be a sequence of analytic sets. 

We first show that $\bigcup_n A_n$ is analytic. For each $n$ there is a
	Polish space $X_n$ and a closed set $F_n\subseteq X\times X_n$ such that $A_n=\Pi_{X}(F_n)$. Let $Y$ be the separated sum of the $X_n$'s. Let $F_\infty=\cup_{n\in\omega} F_n$. 
We claim that $F_\infty$ is closed in $X\times Y$. 
	
	To see this, let ${\la}(x_n, y_n)\ |\ n\in\omega{\ra}$ be a sequence in $X\times Y$ converging to some $(x,y)$ such that  for every $n\in\omega$,
	$(x_n, y_n)\in F_\infty$. Since the $y_n$ converge there is an $N$ such that $y_n\in X_N$ for all large $n$. Hence for all large $n, (x_n, y_n)\in F_N$. Since $F_N$ is closed we must have $(x,y)\in F_N\subseteq F_\infty$.

	Now $x\in\Pi_X(F_\infty)$ if and only if for some $n$,
	$x\in\Pi_X(F_n)$, that is to say if and only if $x\in\cup_{n\in\omega}A_n$. 
	
	For countable intersections, let $V\subseteq\omega^\omega\times(\omega^\omega\times X)$ be a closed universal set for
	$\bpi_{1}^{0}(\omega^\omega\times X)$. Then for every $n\in\omega$ there is an $f_n\in\omega^\omega$ such that 
	\[A_n=\{x\in X\ |\ \exists g\in\omega^\omega\ (f_n,g,x)\in V\}.\]  
	Fix a homeomorphism $g\mapsto \la (g)_n:n\in\omega\ra$ between $\omega^\omega$ and $(\omega^\omega)^\omega$. 
	Let $V^\prime$ be such that 
	\[(g,x)\in V^\prime \Leftrightarrow \ \forall n\ (f_n,(g)_n,x)\in V\]
	Then $V^\prime$ is closed as $\{(g,x)\ |\ (f_n,(g)_n,x)\in V\}$ is closed. 
	
	Now $x\in\cap_{n\in\omega}A_n$ if and only if for every $n$ there is an $h_n$ such that $(f_n,h_n,x)\in V$.
	Suppose that there is such a sequence $\la h_n:n\in\omega\ra$. Let $h$ be such that $(h)_n=h_n$ for every $n$.  Then $(h,x)\in V'$.   On the other hand if $(h,x)\in V'$ then for all $n, (f_n, (h)_n, x)\in V$.  Thus
	$x\in\cap_{n\in\omega}A_n$ if and only if there is a $g$ such that $(g,x)\in V^\prime$, whereby 
	$\cap_{n\in\omega}A_n$ is analytic.
\qed

\begin{lemma}
	Let $X,Y$ be Polish spaces.
	If $A\subseteq X$ is an analytic subset of $X$ and $f:Y\rightarrow X$ is continuous, then $f^{-1}(A)$ is analytic in $Y$. 
\end{lemma}
\pf
	Let $C\subseteq \omega^\omega\times X$ be closed such that $A=\bpi_X(C)$. Define $D\subseteq \omega^\omega\times Y$ to be the set
	$D=\{(w,y)\ |\ (w,f(y))\in C\}$. Notice that $D$ is a closed set and that 
	\begin{eqnarray*}
	y\in f^{-1}(A) & \Leftrightarrow & f(y)\in A\\
	               & \Leftrightarrow & \exists w\ (w,f(y))\in C\\
				   & \Leftrightarrow & \exists w\ (w,y)\in D
	\end{eqnarray*}
so that $f^{-1}(A)$ is the projection of a closed subset of a Polish space.

\qed

It is easy to see that the continuous image of an analytic set is analytic; in particular, the projection of an analytic set is analytic.

The following theorem shows that in a perfect Polish space  there are strictly more analytic sets than Borel sets. 
\begin{theorem}
	If $X$ is a perfect Polish space, then there is an analytic set $A\subseteq X$ such that $X\setminus A$ is not analytic.
\end{theorem}
\pf
	It is enough to show the result for $X=2^{\omega}$. Let $U\subseteq 2^{\omega}\times2^{\omega}$ be universal for $\bsigma_1^1 (2^{\omega})$. 
	Let $f:2^{\omega}\rightarrow2^{\omega}\times2^{\omega}$ be the function $f(x)=(x,x)$. Since $U$ is analytic, $f^{-1}(U)$ is analytic. Let $A=f^{-1}(U)$. 
	Suppose for a contradiction that $D=X\setminus A=\{x\ |\ (x,x)\not\in U\}$ is analytic. Then there is a $y\in2^{\omega}$ such that $D=U_y$. 
Then again:
	\[y\in D \iff (y,y)\in U\iff y\notin D.\]
\qed
 
 \subsection{Understanding analytic sets in terms of trees}

 In Proposition \ref{tree representation of closed sets},  we understood closed subsets of $A^\om$ as the collection of branches through a pruned tree. We now continue this analysis for analytic sets. 
 
 Given sets $A, A'$ we can identify sequences in $(A\times A')^{<\om}$ with pairs of sequences $(\sigma, \tau)$ with $\sigma\in A^{<\om}$ and $\tau\in (A')^{<\om}$. 
Moreover for an arbitrary tree $T\subseteq (A\times A')^{<\om}$, we can identify elements of $[T]$ with pairs of elements in $A^\om\times (A')^{\om}$. For such a tree, let 
 \[p[T]=\{x\in A^\om|(\exists y\in (A')^\om)(x,y)\in [T].\]
By Proposition \ref{laundry list}, if $S\subseteq \omom$ is analytic there is a closed set $C\subseteq \omom\times \omom$ such that $S=\{f:(\exists g)(f,g)\in C\}$. Since $\omom\times\omom$ is canonically homeomorphic to $(\om\times\om)^\om$ we can apply Proposition \ref {tree representation of closed sets} to see that there is a pruned tree $T\subseteq (\om\times\om)^{<\om}$ (with $A=\om\times\om$) such that $C=[T]$.  
Then
 \begin{eqnarray*}
f\in S& \mbox{iff} & (\exists g)(f,g)\in C\\
& \mbox{iff}& (\exists g)(f,g)\in [T]\\
& \mbox{iff} & f\in p[T].
\end{eqnarray*}

\subsection{The CH for analytic sets}

It is a somewhat strange fact that all concrete (e.g. definable) classes of sets satisfy a strong version of CH, at least assuming the appropriate axioms. We present here two early examples of this type of theorem. We first present the theorem for closed sets to illustrate the method and then for analytic sets. The attentive reader will note that Theorem \ref{perfect for closed sets} implies the same result for Borel sets using Theorem \ref{retopologize}.

If $(X, \tau)$ is Polish and $A\subseteq X$ then we say that $A$ \emph{contains a perfect set} if there is a closed (with respect to $\tau$) non-empty perfect subset $P\subseteq A$. This is equivalent to the existence of a continuous, 1-1 map $\phi:2^\om\to X$ with the range $\phi$ included in $A$.

\begin{theorem}\label{perfect for closed sets}
Let $(X,\tau)$ be a Polish space and suppose that $C\subseteq X$ is closed and uncountable. Then $C$ contains a perfect set.
\end{theorem}
\pf Let $d$ be a complete compatible metric. Let $\la U_i:i\in\om\ra$ be a base for $\tau$. We build a scheme of open sets $\la O_s:s\in 2^{<\om}\ra$ such that:
\begin{enumerate}
\item $O_s\supseteq \overline{O}_{s\cat i}$,
\item $O_{s\cat 0}\cap O_{s\cat 1}=\emptyset$,
\item $O_s$ has diameter less than or equal to ${1\over lh(s)}$,
\item $O_s\cap C$ is uncountable.

\end{enumerate}
We do this inductively, setting $O_{\emptyset}=X$. To complete the construction we need the following:
\medskip

\noindent{\bf Claim:}
Given $O_s$ there are $O_{s\cat 0}$ and $O_{s\cat 1}$ that satisfy 2-4.

\pf(Claim) If not we can choose a decreasing sequence of open sets $\la V_n:n\in\om\ra$ such that $V_0=O_s$, $diam(V_n)<{1\over n}, |V_n\cap C|>\omega$ and $(O_{s}\setminus V_n)\cap C$ is countable. But then 
\[O_s\cap C\subseteq\bigcup_{n\in\omega} ((O_s\setminus V_n)\cap C)\cup\bigcap_{n\in\om}V_n.\]
Since $diam(V_n)\to 0$, the intersection $\bigcap_nV_n$ has at most one element. But this implies that $O_s\cap C$ is countable, a contradiction.\qed
Continuing the proof of Theorem \ref{perfect for closed sets}, we first remark that if $a\in 2^\omega$, then $\bigcap_{n}O_{a{\rest}n}$ is a singleton and a member of $C$. To see this, for each $n$ choose $x_n\in O_{a{\rest}n}\cap C$. Then $\la x_n\ra$ is a Cauchy sequence of elements converging to some $y\in C$. For each $n$, $\la x_m:m>n\ra\subseteq O_{a{\rest}n}$, hence $y\in\overline{O}_ {a{\rest}n}$. Since the $O$'s are nested, $y\in \bigcap_{n}O_{a{\rest}n}$.

Define $\phi(a)=a_x$ where $x\in \bigcap_nO_{a{\rest}n}$. Then $\phi$ is clearly a 1-1 continuous map. \qed

The following is also true, but requires a finer analysis.
\begin{exercise} Let $X$ be a Polish space.
Show that if $C\subseteq X$ is closed, then there is a countable set $S$ and a (possibly empty) closed perfect set $P$ disjoint from $S$ such that $C=S\cup P$, and that this decomposition is unique.

(Hint: remove the isolated points of $C$. Then remove the isolated points of the result. Keep going transfinitely.)

\end{exercise}

We now try to understand the  analysis of Theorem \ref {perfect for closed sets} for analytic sets. Suppose that $A\subseteq \omom$ is analytic and uncountable. Then there is a closed set $C\subseteq \omom\times\omom$ such that $A=\Pi_X(C)$, where $\Pi_X$ is the projection to the first coordinate. Since $A$ is uncountable, $C$ must be uncountable and hence $C$ contains a perfect set. Unfortunately this does not imply that $A$ contains a perfect set. (For example  $C=\{a\}\times\omom$ contains a perfect set, but its projection is $\{a\}$.) Instead we need to do the previous argument more carefully.

If $T\subseteq (S)^{<\om}$ and $\sigma\in (S)^{<\om}$ we let 
\[T_\sigma=\{\sigma'\in T: \sigma\subseteq \sigma'\mbox{ or }\sigma'\subseteq \sigma\}.\]

\begin{theorem}\label{perfect set theorem for analytic sets}
Suppose that $A\subseteq \omom$ is an uncountable analytic set. Then $A$ contains a perfect subset.
\end{theorem}
\pf Let $C\subseteq \omom\times \omom$ be a closed set such that $A=\Pi_X(C)$. Let $T\subseteq (\om\times\om)^{<\om}$ be a tree such that $A=p[T]$.

\noindent{\bf Claim} Suppose that $p[T_\sigma]$ is uncountable. Then there are $\sigma_0, \sigma_1\in T_\sigma$ such that
\begin{enumerate}
\item $p[T_{\sigma_i}]$ is uncountable for each $i$ and 
\item $p[T_{\sigma_0}]\cap p[T_{\sigma_1}]=\emptyset$.
\end{enumerate}
Suppose that $\sigma=(\rho, \tau)$. A simple counting argument shows that there must be $\eta\supseteq \rho$ and $m\ne n$
$|p[T_\sigma]\cap [\eta\cat m]|>\omega$ and $|p[T_\sigma]\cap [\eta\cat n]|>\omega$. Another counting argument shows that there must be $\tau_0\supseteq \tau$ and $\tau_1\supseteq \tau$ such that 
both $p[T_{\eta\cat m,\tau_0)}]$ and $p[T_{(\eta\cat n, \tau_1)}]$ are uncountable.
Letting $\sigma_0=(\eta\cat m, \tau_0)$ and $\sigma_1=(\eta\cat n, \tau_1)$ we see the claim.

We can now use the claim to produce a collection $\la \sigma_s:s\in 2^{<\om}\ra$ such that:
\begin{enumerate}
\item $\sigma_{s\cat i}\supseteq \sigma_s$,
\item $p[T_{\sigma_s}]$ is uncountable for all $s$  and 
\item $p[T_{\sigma_{s\cat0}}]\cap p[T_{\sigma_{s\cat 1}}]=\emptyset$.

\end{enumerate}
Define $\phi:2^\om\to [T]$ by $\phi(a)\rest n=\sigma_{a\rest n}$. Then $\phi$ is a continuous map and for $a\ne b\in 2^\om$, we see that $\Pi_X(\phi(a))\ne \Pi_X(\phi(b))$ as needed.\qed
\begin{exercise}
Show that if $X$ is a Polish space and $A\subseteq X$ is an uncountable analytic set then $A$ contains a perfect subset.
\end{exercise}

The next exercise extends \ref{borel bijections one}:
\begin{exercise}\label{borel bijections two}
Suppose that $X$ is an uncountable Polish space. Then there is a bijection  $f:X\to \omom$ such that $f$ and $f^{-1}$ are Borel measurable.
\end{exercise}

\subsection{Operation $\mathcal A$}

This section is about Suslin's operation $\mathcal A$. 

\begin{definition} Let $X$ be perfect Polish and $\mathcal S=\la C_s:s\in \omega^{<\omega}\ra$ be a scheme. 
Define $\mathcal A(\mathcal S)=\bigcup_{f\in\omom}\bigcap_{n\in\omega}C_{f\rest n}.$ If $\mathcal R\subseteq P(X)$  is any collection of sets  we let $\mathcal A(\mathcal R)$ be the result of applying operation $\mathcal A$ to schemes of sets in $\mathcal R$.

\end{definition}
Note the form of the operation:
\[x\in \mathcal A(\mathcal S)\mbox{ iff }(\exists f\in \omom) (\forall n\in \omega) x\in C_{f{\rest}n}.\]

\begin{theorem}\label{operation A is idempotent}
Let $X$ be Polish and $\mathcal R\subset P(X)$. Then $\mathcal A(\mathcal A(\mathcal R))=\mathcal A(\mathcal R)$.
\end{theorem}
\pf It suffices to show that $\mathcal A(\mathcal A(\mathcal R))\subseteq \mathcal A(\mathcal R)$.


Let $\la B_\sigma:\sigma\in \omega^{<\om}\ra$ be a sequence of $\mathcal A(\mathcal R)$ sets. We need to see that $A=\bigcup_{f\in \omom}\bigcap_{n\in\om}B_{f{\rest}n}\in \mathcal A(\mathcal R)$. Write each $B_\sigma=\bigcup_{g\in\omom}\bigcap_{m\in\om}C_{\sigma, g{\rest}m}$.

Note that 
\begin{eqnarray}
A&=&\{x:(\exists f\in\omom)(\forall n)(\exists g\in \omom)(\forall m)(x\in C_{f{\rest}n, g{\rest}m})\}\\
 &=&\{x:(\exists f\in\omom)(\exists \la g_n:n\in\om\ra)(\forall n)(\forall m)(x\in  C_{f{\rest}n, g_n{\rest}m})\}\label{4.2}\\
 &=&\{x:(\exists h)(\forall n)(\forall m)(x\in  C_{h^*{\rest}n, (h)_n\rest m})\}
\end{eqnarray}

where $h\mapsto (h^*, \la (h)_n:n\in\om\ra)$ is a homeomorphism of $\omom$ with $\omom\times (\omom)^\om$. Unfortunately this is not \emph{exactly} the form we need.

Let $\la \cdot, \cdot \ra:\om\times \om\to \om$ be a bijection such that for all $m, n$, we have $n\le \la n, m\ra, m\le\la n,m\ra$ and $\la m, \cdot\ra:\om\to \om$ is order preserving. Let $n\mapsto ((n)_0,(n)_1)$ be the inverse map. 

We define a special homeomorphism from $\phi:\omom\times (\omom)^\om\to \omom$ by
sending $\phi(f,\la g_n:n\in\om\ra)=h$ where 
\[h(k)=\la f(k), g_{(k)_0}((k)_1)\ra.\]
We note that if we know $h\rest \la n,m\ra$ then we can recover $f\rest n$ and $g_n\rest m$.
Given $\tau\in\om^{<\om}$ of length $k=\la n, m\ra$ we treat it as $h\rest\la n, m\ra$ and define 
\[D_\tau=C_{f\rest n, g_n\rest m}.\]
We now check that $\mathcal A(\la D_\tau:\tau\in \om^{<\om}\ra)=A$. Suppose that $x\in A$. Let $f$ and $\la g_n:n\in\om\ra$ be as in Equation \ref{4.2}. Let $h=\phi(f, \la g_n:n\in\om\ra)$. Then 
\[
(\forall n)(\forall m)(x\in C_{f\rest n, g_n\rest m})\]
so
\[(\forall k)(x\in D_{h\rest k}).\]
On the other hand, suppose that $x\in \mathcal A(\la D_\tau\ra)$. Then there is an $h$ such that for all $k, x\in D_{h\rest k}$. Let $(f, \la g_n:n\in\om\ra)=\phi^{-1}(h)$. Then:
\[(\forall n)(\forall m)(x\in C_{f\rest n, g_n\rest m}).\]\qed

\begin{theorem}\label{more a}
Let $X$ be perfect Polish and $A\subseteq X$. Then the following are equivalent:

\begin{enumerate}
\item $A$ is analytic
\item $A$ is the result of applying operation $\mathcal A$ to a scheme of Borel sets
\item $A$ is the result of applying operation $\mathcal A$ to a scheme of closed sets.
\end{enumerate}
In particular, the analytic sets are closed under operation $\mathcal A$.
\end{theorem}
\pf In view of Theorem \ref{operation A is idempotent}, it suffices to show that 1 is equivalent to 3. Suppose $A$ is analytic. Then there is a closed set $C\subseteq X\times\omom$ such that $A=\Pi_X(C)$. Let $C_s$ be the closure of $\{x:(\exists y)((x,y)\in C\mbox{ and } y\in [s])\}$. We claim that $A=\mathcal A(\la C_s:s\in \om^{<\om}\ra)$. Clearly $A\subseteq \mathcal A(\la C_s:s\in \om^{<\om}\ra)$. So suppose that $x\in \mathcal A(\la C_s:s\in \om^{<\om}\ra)$. Then we can find a $g\in \omom$ such that for all $n$,  there is a $y_n, (x, y_n)\in \overline{\{(x,y):(x, y)\in C\mbox{ and }y\rest n=g\rest n\}}$. Hence we can find $x_n, y_n$ such that $x_n\to x$ and $(x_n, y_n)\in C$ and $y_n\rest n=g\rest n$.  Then $y_n\to g$, so $\la x_n, y_n\ra\to (x, g)$. Since $C$ is closed, we have $(x,g)\in C$ and $x\in A$.

Now suppose that $A=\mathcal A(\la D_s:s\in \om^{<\om}\ra)$ where each $D_s$ is closed. Define $C=\{(x, f):(\forall n)(x\in D_{f\rest n})\}$. Then $C$ is closed and $A=\Pi_X(C)$.\qed

Inspecting the proof of Theorem \ref{more a}, we see that if we let $B_s=\{x:(\exists y)((x,y)\in C\mbox{ and } y\in [s])\}$ then $A=B_{\emptyset}=\mathcal A(\la B_s:s\in\om^{<\om})=\mathcal A(\la C_s:s\in\om^{<\om}\ra)$. This allows us to give an easy proof of:
\begin{theorem}\label{analytic sets are measurable}
Let $\mu$ be a regular, non-atomic probability measure on $X$ and $ A$ be an analytic set. Then $A$ is $\mu$-measurable.
\end{theorem}

\pf We use the notation of the proof of Theorem \ref{more a}.  As $\mu$ is regular we can find Borel sets $D_s\supseteq B_s$ such that the outer measure of $B_s$ is $\mu(D_s)$. One way of saying this is that if $E\subset D_s\setminus B_s$ is measurable, then $\mu(E)=0$. Without loss of generality we can assume that $B_s\subseteq D_s\subseteq C_s$ and thus $A=B_\emptyset=\mathcal A(\la D_s:s\in \om^{<\om}\ra)$.

We claim that $\mu(D_\emptyset\setminus B_\emptyset))=0$. This clearly suffices. To see that this is true, note that 
\[D_\emptyset\setminus B_\emptyset\subseteq \bigcup_{s\in\om^{<\om}}(D_s\setminus \bigcup_nD_{s\cat n}).\]
This follows since if $x\in D_{\emptyset}$ but $x\notin  \bigcup_{s\in\om^{<\om}}(D_s\setminus \bigcup_nD_{s\cat n})$ we can inductively build $f\in \omom$ such that $x\in \bigcap_n D_{f\rest n}$. However, we know that $\mathcal A(\la D_s:s\in \om^{<\om}\ra)=\mathcal A(\la B_s:s\in \om^{<\om}\ra)$. Hence $x\in B_\emptyset$.

Thus it suffices to show that for each fixed $s\in \om^{<\om}, \mu(D_s\setminus \bigcup_nD_{s\cat n})=0$. Let $E\subseteq D_s\setminus \bigcup_nD_{s\cat n}$ be Borel. Then 
\[E\subseteq D_s\setminus \bigcup_nD_{s\cat n}\subseteq D_s\setminus \bigcup_n B_{s\cat n}\subseteq D_s\setminus B_s.\]
Hence $\mu(E)=0$.\qed

\begin{exercise}
Show that if $A\subseteq X$ is analytic and $X$ is Polish, then $A$ has the Property of Baire.
\end{exercise}

\subsection{Lusin's separation theorem and $\bdelta^1_1$ sets}

\begin{theorem}[Lusin's Separation Theorem]
Let $X$ be a Polish space, $A$ and $B$ be disjoint analytic subsets of $X$. Then there is a Borel set $C$ with $A\subseteq C$ and $C\cap B=\emptyset$.

\end{theorem}

\pf
	Let $f,g$ be continuous functions from $\omega^\omega$ to $X$ such that $A=range(f)$ and $B=range(g)$. For $s\in\omega^{<\omega}$, 
	let $A_s=f([s])$ and $B_s=g([s])$. Then $A_s=\bigcup_{i\in\omega}A_{s\cat i}$ and
	$B_s=\bigcup_{j\in\omega}B_{s{\cat}j}$. 
\bigskip

\noindent{\bf Observation:}
	Suppose that for all $i,j\in\omega$ there is a Borel set $C_{ij}$ such that $A_{\sigma_{{\cat}}i}\subseteq C_{ij}$ and $B_{\tau{\cat}j}
	\cap C_{ij}=\emptyset$. Then $A_{\sigma}\subseteq \bigcup_i\bigcap_j C_{ij}$ and $\bigcup_i\bigcap_j C_{ij}\cap B_{\tau}=\emptyset$.

\bigskip
	
	Let us assume, for a contradiction, that $A$ and $B$ cannot be separated by any Borel set $C$. Using the claim, we can inductively build sequences ${\la}\sigma_k\ |\ k\in\omega{\ra}$ and ${\la}\tau_k\ |\ k\in\omega{\ra}$ with the following properties:
	\begin{enumerate}
	\item $\sigma_k,\tau_k\in\omega^k$
	\item $\sigma_k\subseteq\sigma_{k+1}\ \tau_k\subseteq\tau_{k+1}$
	\item $A_{\sigma_k} \mbox{ can't be separated from $B_{\tau_k}$ by any Borel set}$
	\end{enumerate}

	Let $x=\bigcup_{k\in\omega}\sigma_k$ and $y=\bigcup_{k\in\omega}\tau_k$.
	Consider $f(x)\in A$ and $g(y)\in B$. We can find open sets $U,V$ in $X$ with $f(x)\in U$, $g(y)\in V$ and $U\cap V=\emptyset$. 
	
	By the continuity of $f$ and $g$, there is a $k\in\omega$ with 
	\[f([\sigma_k])\subseteq U\mbox{ and }g([\tau_k])\subseteq V\]
	Then $A_{\sigma_k}$ and $B_{\tau_k}$ can be separated by these $U,V$, both of which are Borel, a contradiction.
\qed

\begin{corollary}[Suslin's Theorem]\label{suslin theorem} Let $A\subseteq X$. 
	If   both $A$ and $X\setminus A$ are analytic, then $A$ is Borel.
\end{corollary}
\pf
	Use the Lusin separation theorem on $A$ and $X\setminus A$.
\qed

\begin{corollary}\label{analytic functions are Borel}
	Suppose that $X,Y$ are Polish spaces, and that $f:X\rightarrow Y$. Then the following are equivalent:
	
	\begin{enumerate}
	\item $f$ is Borel measurable.	
	\item $graph(f)$ is a Borel subset of $X\times Y$.
	\item $graph(f)$ is an Analytic subset of $X\times Y$.
	\end{enumerate}
	
\end{corollary}
\pf
	$(1)\Rightarrow (2)$ is from Proposition (\ref{graphborel}). $(2)\Rightarrow (3)$ is trivial. It remains to show that $(3)\Rightarrow (1)$.
	Let $O\subseteq Y$ be open. We will show that $f^{-1}(O)$ is Borel in $X$. In particular, we show that it and its complement are both analytic.
	Now $f^{-1}(O)=\{x\in X\ |\ \exists y\ (y=f(x)\land y\in O)\}$ and $X\setminus f^{-1}(O)=\{x\in X\ |\ \exists y\ (y=f(x)\land y\not\in O)\}$, both of which are analytic sets.
\qed

\begin{theorem}
	Let $X$, $Y$ be Polish spaces and $B\subseteq X$ be a Borel set. If $f:X\rightarrow Y$ is continuous and $f\rest B$ is injective, then
	$f(B)\subseteq Y$ is Borel.
	
\end{theorem}
\pf	By Theorem \ref {retopologize} , we can retopologize $X$ to make $B$ a closed set in a Polish topology finer than the original topology on $X$.  By Theorem \ref{cont bij}, there is a closed set $C\subseteq \omom$ and a function $g:\omom\to B$ that is continuous with respect to the finer topology such that $f\rest C$ is a one to one map onto $B$. Since $g$ is continuous with respect to the finer topology, $g$ is continuous with respect to the original topology on $X$. Composing $f$ with $g$ we get a continuous map of $\omom$ to $Y$ and a closed set $C$ such that $f\circ g$ is a one one map from $C$ to $f(B)$.
Thus	 it is enough to consider the case where $X=\omega^\omega$ and $B$ closed. Let $T\subseteq \om^{<\om}$ be a pruned tree such that $B=[T]$.

	Let $A=f(B)$, and for every $s\in T$, let $A_s=f(B\cap[s])$. Fix a complete compatible metric on $Y$.
	 Then
\begin{enumerate}
\item $A_{s\cat i}\cap A_{s\cat j}=\emptyset$ for $i\ne j$ since $f$ is 1-1,
\item $A_{s\cat i}\subseteq A_s$,
\item each $A_s$ is analytic,
\item for every $x\in B$, diam$(A_{x\rest n})\to 0$.
\end{enumerate}
 Using the Lusin separation theorem to separate $A_{s\cat i}$ from $\bigcup_{j\ne i}A_{s\cat j}$ we can recursively construct a scheme of Borel sets $\la A'_s:s\in\om^{<\om}\ra$ such that 
\begin{enumerate}
	\item[a.] $A_s\subseteq A^{\prime}_s$.
	\item[b.] $A^{\prime}_{s{\cat} i}\cap A^{\prime}_{s{\cat} j}=\emptyset$
	\item[c.] $A^{\prime}_{s{\cat} i}\subseteq A^{\prime}_{s}$
	\end{enumerate}

	
If we set $A^*_s=A'_s\cap \overline{A_s}$ then we get another scheme of Borel sets that satisfies (a)-(c) and for all $x\in B$, diam$(A^*_{x\rest n})\to 0$. The following claim suffices to show the theorem.
	\begin{claim}
		\[f(B)=\bigcap_{k\in\omega}\bigcup_{s\in(T\cap \omega^k)} A^{\ast}_s\]
	\end{claim}

	\pf[Claim]
		Let $x\in B$. Then $f(x)\in \bigcap_{k\in\omega} A_{x{\rest} k}\subseteq \bigcap_{k\in\om}\bigcup_{s\in(T\cap\om^k)}A^*_s$. Hence $f(B)\subseteq  \bigcap_{k\in\om}\bigcup_{s\in(T\cap\om^k)}A^*_s$.
				
		Now assume that $y\in \bigcap_{k\in\omega}\bigcup_{s\in(T\cap\omega^k)} A^{\ast}_s$. 
		Then for each $k$ there is a unique  $s_k\in\om^k$ with $y\in A^\ast_{s_k}$. By our tree construction,
		$s_{k}\subseteq s_{k+1}$, so we can let $x=\bigcup_{k\in\omega}s_k$. 
		
Since $x$ is a branch through $T$, $x\in B$. We claim that $f(x)=y$. If not, let $z=f(x)$. Then $z\in \bigcap_{k\in \om}A^*_{x\rest k}$. Since diam$(A^*_{x\rest k})\to 0$, we have $d(y,z)=0$. 
%
\qed
We introduce the following definition here, but will only use it later.

\begin{definition}
Let $X$ be a Polish space, and $A\subseteq X$. $A$ is called co-analytic if $X\setminus A$ is analytic. The class of co-analytic sets of 
$X$ is denoted by $\bpi_1^1(X)$. We denote $\bpi^1_1\cap \bsigma^1_1$ by $\bdelta^1_1$. 
\end{definition}

In this notation: 
\[ \bsigma_1^1(X)\cap\bpi_1^1(X)=\bdelta^1_1(X)=\mbox{Borel}(X).\]

\begin{remark}
In light of the results in this and earlier sections the properties of analytic and co-analytic depend largely on the Borel structure rather than the particular topology on a space $X$. Thus we can rephrase many of these results in terms of \emph{standard Borel spaces}.
\end{remark}

\section{Reductions}

\begin{definition}
	Let $X,Y$ be Polish spaces with $A\subseteq X$ and $B\subseteq Y$. A function $f:X\rightarrow Y$ is called a reduction if and only if
	\[\forall x\in X\ (x\in A\Leftrightarrow f(x)\in B)\]
	We write $A\leq B$. 
	
A reduction $f$ is called a Borel reduction provided $f$ is Borel measurable and write $A\leq_{\mathbb{B}} B$.
	If $f$ is  continuous,
	$f$ is called a continuous reduction 	and we write $\le_c$.
\end{definition}
We note that ``reductions" reduce the question of membership in $A$ to membership in $B$ and that the relations $\le_{\mathbb B}$ and $\le_c$ are both transitive.
\begin{definition}
	Let $X$, $Y$ be Polish spaces, $\Gamma$ be a collection of subsets of $X$, and $B$ be a subset of $Y$. 
	Then $B$ is \emph{complete} for $\Gamma$ via Borel (or continuous) reduction if and only if for all $A\in\Gamma$ there is a Borel (or continuous) reduction of $A$ to $B$. 
	
If $\Gamma$ is a collection of subsets of Polish spaces, 	we will say that $B$ is a \emph{complete $\Gamma$ set} iff $B$ is a $\Gamma$ set and is complete for the $\Gamma$ subsets of each Polish space. 
	
\end{definition}
\begin{example}
	\begin{enumerate}
	\item Let $U\in \bsigma_\alpha^0(\omega^\omega\times X)$ be universal for $\bsigma_\alpha^0(X)$. Then $U$ is complete for $\bsigma_\alpha^0$ subsets of $X$ via continuous reduction.
	\item Let $X=\omega^\omega$. Then there is a $\bsigma^0_\alpha$ set $U'\subseteq \omega^\omega$ that is complete via continuous reductions for 
	the $\bsigma^0_\alpha$ subsets of $\omega^\omega$. 
	\item The same results are true for $\bpi^0_\alpha$, $\bsigma^1_1$ and $\bpi^1_1$.
	\end{enumerate}
\end{example}
\pf To see the  first claim in the example, we first take a $\bsigma^0_\alpha$ set $A\subseteq X$. Then there is an $x_0\in\omom$ such that $A=U_{x_0}$. Define a function $R:X\to \omom\times X$ by $y\mapsto (x_0, y)$. It is easy to check that $R$ is a reduction and continuous.

For the second claim, let $F:\omom\times \omom \to \omom$ be a homeomorphism.  Let $U'=F[U]$. Given $A\subseteq\omom$ be $\bsigma^0_\alpha$ and suppose that $A=U_{x_0}$. Define $R:\omom\to \omom$ by $R(y)=F(x_0, y)$. As before $R$ is a continuous reduction.

The third claim is proved along the same lines by taking $U$ to be universal for different pointclasses.\qed

\subsection{Reducing analytic sets to ill-founded trees}

\begin{definition}
Let $A$ be a set and $T\subseteq A^{<\omega}$ be a tree. We call $T$ \emph{well-founded} if and only if $T$ has no infinite branches. If $T$ does have
an infinite branch, then we call $T$ \emph{ill-founded}.
\end{definition}

If $T\subseteq A^{<\om}$ is a tree, we define $<_T$ on $T$ by setting $\sigma<_T\tau$ iff $\sigma\supset \tau$. With this partial ordering $T$ is well-founded iff the ordering $<_T$ is well-founded in the sense of Definition \ref {well founded relations}. As a consequence $T$ is well-founded iff there is a function:
\[\rho:T\to OR\]
such that $\sigma\subset \tau$  and $\sigma\ne \tau$ implies $\rho(\sigma)>\rho(\tau)$.
We call the least ordinal $\alpha$ such that there is a $<_T$-preserving function $\rho:T\to \alpha$ the \emph{height} of $T$. We note that if $T\subseteq A^{<\om}$ then the height of $T$ has cardinality at most $|A|+\aleph_0$.

\begin{exercise}
Let $\alpha$ be a countable ordinal. Then there is a well-founded tree $T\subseteq \om^{<\om}$ such that $T$ has height $\alpha$. 
\end{exercise}

Now, if $A\subseteq \om^\omega$ is analytic, then there is a closed set $C\subseteq \omega^\omega\times\omom$ such that 
$A=\{x\ |\ (\exists f\in\omega^\omega)\ (f,x)\in C\}$. Let $T\subseteq \omega^{<\omega}\times\om^{<\omega}$ be defined by:
\[T=\{(s,t)\ |\ lh(s)=lh(t)\ \land\ ([s]\times[t])\cap C\not=\emptyset\}\]
Then $C=\{(f,x)\ |\ \forall n (f \rest  n,x \rest  n)\in T\}$

Given $x\in \omom$, let $ {T^{x}}=\{s\ | (s,x \rest lh(s))\in T\}$. Then $ {T^{x}}$ is a tree and 
\begin{eqnarray*}
	x\in A & \Leftrightarrow & \exists f\forall n\ (f \rest  n,x \rest  n)\in T\\
	       & \Leftrightarrow & \exists f\forall n\ f \rest  n\in  {T^{x}}
  \end{eqnarray*}
Thus $x\in A$ if and only if $ {T^{x}}$ is ill-founded.

\medskip
As we saw in Example \ref {space of infinite trees} , the collection of trees can be represented as a closed subset of $2^\om$, via an enumeration $\la \sigma_n:n\in \omega\ra$ of $\om^{<\om}$. Hence we can view the collection of trees (or the collection of infinite trees) as a Polish space, which we named $\trees$. For convenience we will view $2^\om$ as $2^{\omega^{<\omega}}$ 

\begin{lemma}
The map $R: \omom \rightarrow\trees$ given by $R(x)= {T^{x}}$ is continuous.
\end{lemma}

\pf
Basic open sets for $\trees$ have the form $B(\sigma_1,\dots,\sigma_k,\tau_1,\dots,\tau_l)=\{T\ | \forall i<k,j<l\ \sigma_i\in T\ \land\ 
 \tau_j\not\in T\}$. 
 
 Fix $\sigma_1,\dots,\sigma_k,\tau_1,\dots,\tau_l$.  Let $n
 =\max_{i<k,j<l} \{lh(\sigma_i),lh(\tau_j)\}$. Suppose that $x, y\in 2^\om$ are such that $x\rest n=y\rest n$. Then $R(x)\in B(\sigma_1,\dots,\sigma_k,\tau_1,\dots,\tau_l)$ iff
 $R(y)\in B(\sigma_1,\dots,\sigma_k,\tau_1,\dots,\tau_l)$. Thus $R^{-1}(B(\sigma_1,\dots,\sigma_k,\tau_1,\dots,\tau_l)$ is a clopen set.
\qed

Hence we have shown that if $A$ is analytic, then there is a continuous function $R:\omega^\omega\rightarrow \trees$ such that $x\in A$ if and 
only if $ {T^{x}}$ is ill-founded.

\begin{theorem}
The set $\{\mbox{ill-founded trees}\}\subseteq \trees$ is an  analytic set  and every analytic subset of $\om^\om$ can be continuously reduced to it.  Moreover if $A\subset X$ is an analytic subset of a Polish space $X$, then there is a Borel reduction of  $A$ to $\trees$.
\end{theorem}

\noindent\emph{Comment:} Kechris has shown that every set that is complete via Borel reductions is also complete via continuous reductions, but we don't prove or use that fact in these notes.
\pf
We have shown that every analytic subset of $\omom$ can be continuously reduced to $\trees$ so what remains is to reduce arbitrary analytic sets to $\trees$.  We can assume that $X$ is uncountable, since every subset of a countable Polish space is an $\mathcal F_\sigma$-set.
\medskip

\noindent{\bf Claim:} 	Let $A\subseteq X$ be an analytic subset of an uncountable Polish space. Let ${\la}U_n\ |\ n\in\omega{\ra}$ be an enumeration of a basis for 
	the topology of $X$. Define a map $c:X\rightarrow 2^\omega$ by $c(x)=\chi_{\{n | x\in U_n\}}$. Then $c$ is a Borel measurable injection.

To see the claim, we note that If $[s]$ is a basic open set in $2^\omega$, then 
	\[c^{-1}([s])=\bigcap_{s(n)=1} U_n\cap\bigcap_{s(m)=0}U_m^c\]
	which is Borel.

	Let $A$, $X$ and $c$ be as in the above proposition. Let $A^\prime=c(A)$, which is an analytic subset of $2^\omega$. Let $R:A' \to \trees$
	be a continuous reduction of $A^\prime$ to ill-founded trees. Then 
	\begin{eqnarray*}
	x\in A & \Leftrightarrow & c(x)\in A^\prime\\
		   & \Leftrightarrow & R(c(x))\mbox{ is ill-founded}
	\end{eqnarray*}
	so that $R\circ c:X\rightarrow \trees$ is a Borel measurable reduction.\qed
\begin{corollary}The set $\{\mbox{well-founded trees}\}\subseteq \trees$ is a  complete co-analytic set.
\end{corollary}
\subsection{Linear Orderings and Well Orderings}

\begin{definition}
	Consider relations $R\subseteq \omega\times\omega$ that are linear orderings. We call $\mathcal{L}\mathcal{O}=\{\chi_R\in2^{\omega\times\omega}\
	|\ R\mbox{ is a linear ordering}\}$. 
	
	Similarly, we write $\mathcal{W}\mathcal{O}$ for the set of characteristic functions of well orderings on $\omega$.
\end{definition}
If $T$ is a tree we want to be able to systematically linearize $<_T$ in such a way that $T$ is well-founded iff its linearization is a well-ordering.  The \emph{Kleene-Brower} ordering does this with a local definition. (In other contexts this is sometimes called the \emph{Lusin-Sierpinski} ordering.)

\begin{definition}
	The $KB$ relation on $\omega^{<\omega}$ is defined as follows:
	
	$s\prec_{KB} t$ if and only if either $t\subseteq s$ or if there is an $n<lh(s)$ such that for every $i<n$ we have that $s(i)=t(i)$ and 
	$s(n)<t(n)$. 
\end{definition}

We can extend the $KB$ relation to a relation on $(\omega\times \omega)^{<\om}$ (or $(\om^n)^{<\om}$) by setting 
\[(\sigma,\tau)\prec_{KB}(\sigma', \tau')\]
\begin{center}
iff
\end{center}
\[(\sigma(0),\tau(0),\sigma(1), \tau(1), \dots \sigma(n-1), \tau(n-1))\prec_{KB}(\sigma'(0),\tau'(0), \dots \sigma'(m-1), \tau'(m-1))\]
where $n=lh(\sigma)$ and $m=lh(\sigma')$.

\begin{proposition}
Let $\alpha$ be any ordinal. Then:

\begin{itemize}
\item $\prec_{KB}$ is a linear ordering of $\alpha^{<\om}$ with greatest element $\emptyset$,
\item  If $T\subseteq \alpha^{<\omega}$ is a tree, then $T$ is well-founded if and only if $(T,\prec_{KB})$ is a well ordering.
\end{itemize}

\end{proposition}
\pf
	To check that $\prec_{KB}$ is a linear ordering is a routine exercise left for the reader. 
	
	For the second clause, suppose that $T$ is ill-founded. Let $f:\omega\rightarrow\omega$ be such that for every $n$, $f{\rest} n\in T$.
	Then $\{f{\rest} n\ |\ n\in\omega\}$ is $\prec_{KB}$-decreasing, so that $\prec_{KB}$ cannot be a well-ordering.
	
	Now suppose that $(T,\prec_{KB})$ is not a well ordering. 
	Let ${\la}\sigma_k\ |\ k\in\omega{\ra}$ be a $\prec_{KB}$-decreasing sequence. Define the function $f:\omega\rightarrow\alpha$ by induction:
	\begin{eqnarray*}
	f(0) & = & \mbox{the least }\beta\mbox{ such that }\exists k\ \sigma_k(0)=\beta\\
	f(n) & = & \left\{\begin{array}{ll}
					\mbox{the least $\beta$} & \mbox{such that $\exists k\ \sigma_{k}{{\rest} n}=f{\rest} n\ \land \sigma_k(n)=\beta$}\\
					\mbox{undefined} & \mbox{otherwise}
					\end{array}\right.
	\end{eqnarray*}
	We show by induction on $n$ that $f{\rest} n$ is the $\prec_{KB}$-least element of 
	\[\{\tau\ |\ \mbox{for some $k$, $\tau\subseteq \sigma_k$ and $lh(\tau)=n$}\}\]
	This will imply that $f(n)$ is always defined.
	
	Suppose we have shown the statement to be true for $n$. Suppose that $f\rest n=\sigma_k\rest n$.  Since $\sigma_{k+1}\prec_{KB}\sigma_k$ and $\sigma_{k+1}\rest n\not\prec_{KB}\sigma_k\rest n$, we must have $\sigma_{k+1}$ strictly extends $\sigma_k$ and hence we know that $f(n)$ exists. Suppose for a contradiction that 
	$\sigma_l{\rest} n+1\prec_{KB} f{\rest} n+1$. Then either $\sigma_l{\rest} n\prec_{KB} f{\rest} n$ (which is
	impossible by our induction hypothesis), or $\sigma_l(n)<f(n)$, which cannot be by the definition of $f$. 
	
	Now from the definition of $f$, it is now clear that $f{\rest} n\in T$ for each $n$. Hence $f$ is an infinite branch through
	$T$.
\qed

For the following result, let ${\la}\sigma_n\ |\ n\in\omega{\ra}$ be an enumeration of $\omega^{<\omega}$ such that $\sigma_n\not\subseteq\sigma_m$
whenever $m<n$.
Define $R:\trees\rightarrow \mathcal{L}\mathcal{O}$ to be:
\[ R(T)(m,n)=\left\{\begin{array}{ll}
					1 & \mbox{if $\sigma_n,\sigma_m\in T$ and $\sigma_n\prec_{KB}\sigma_m$}\\
					0 & \mbox{otherwise}
					\end{array}\right.\]
\begin{theorem}
$R$ is a continuous reduction from  $\{$well-founded trees$\}\subseteq \trees$ to $\mathcal{WO\subseteq LO}$. \end{theorem}
\pf
	That $R$ is a reduction is clear from the previous proposition.
	We need to see that $R$ is continuous.

Let $s:r\times r\to 2$ determine a basic open interval in $\mathcal{LO}$. Then
$R^{-1}(\la s\ra)$ is:
\begin{equation*}\{T:\mbox{for all }m,n<r((\sigma_m, \sigma_n\in T \mbox{ and } \sigma_m\prec_{KB}\sigma_n)
\mbox{ iff }
s(m,n)=1)\}
\end{equation*}	

Since $\prec_{KB}$ is independent of $T$, if $s$ coincides with $\prec_{KB}$ then $R^{-1}{[s]}=\{T:n=0\mbox{ or } (\mbox{for all }n<r)\sigma_n\in T \mbox{ iff } s(n,0)=1\}$.
\qed
\begin{corollary}
The set $\mathcal{W}\mathcal{O}$ is a co-analytic set that is complete via continuous reduction for co-analytic subsets of $\omega^\omega$.
\end{corollary}
\subsection{Rank on the family of  trees}

\begin{definition}

Let $T\subseteq \omega^{<\omega}$ be a tree.
\begin{enumerate}
	\item Define the following sequence of trees ${\la}T_\alpha\ |\ \alpha\in\omega_1{\ra}$:
\begin{eqnarray*}
	T_0 & = & T\\
	T_{\alpha+1} & = & \{\sigma\ |\ \sigma \mbox{ is not a terminal node of $T_\alpha$}\}\\
	T_\lambda & = & \bigcap_{\beta<\lambda} T_\beta\mbox{\ \  if $\lambda$ is a limit ordinal}
\end{eqnarray*}
	\item Define the following function $\rho_T: \omega^{<\omega}\rightarrow \omega_1\cup\{\infty\}\cup\{-1\}$ by:
	\[\rho_T(\sigma)=\left\{\begin{array}{ll}
							-1 & \mbox{if $\sigma\not\in T$}\\
							\mbox{the least $\alpha$} & \mbox{such that $\sigma\not\in T_{\alpha+1}$}\\
							\infty & \mbox{if no such $\alpha$ exists}
						\end{array}\right.\]

\item if $T$ is well-founded and non-empty, we let $\rho(T)=\rho_T(\emptyset)$ and $-1$ otherwise.
	\item If $T$ is ill-founded we let $\rho(T)=\infty$.
	\item For any given $\sigma\in\omega^{<\omega}$, let $T(\sigma)=\{\tau\ |\ \sigma{\cat}\tau \in T\}$.
\end{enumerate}
\end{definition}

Notice that since $\abs{T}=\aleph_0$, there is an $\alpha<\omega_1$ such that for all $\beta>\alpha$, $T_\beta=T_\alpha$. We sometimes denote this $T_\alpha$ by $T_\infty$.
\begin{exercise}
\
\begin{enumerate}

\item \[\rho_T(\sigma)=\left\{\begin{array}{ll}
								\sup_{i\in\omega}\{\rho_T(\sigma{\cat} i)+1\ |\ \sigma{\cat} i\in T\} & \mbox{if $\sigma\in T$}\\
								-1 & \mbox{otherwise}
							\end{array}\right.\]
\item If $\sigma\subseteq\tau$, then $\rho_T(\sigma)\geq\rho_T(\tau)$.
\item If $\sigma\subset\tau$, $\sigma\ne \tau$ and $0\leq\rho_T(\tau)<\infty$, then $\rho_T(\sigma)>\rho_T(\tau)$.
\item Show that for every $\sigma\in\omega^{<\omega}$, the function $\phi_\sigma: {{\trees}}\rightarrow {{\trees}}$ given by 
\[\phi_\sigma(T)=T(\sigma)\]
is continuous and
 $\rho(T(\sigma))=\rho_T(\sigma)$.
\end{enumerate}
\end{exercise}

\begin{definition}
\begin{enumerate}
\item We denote by $\mathcal{W}\mathcal{F}_\alpha$ the set of trees $T$ with $\rho(T)<\alpha$. We say these are well founded trees of rank $\alpha$.
\item Let $S,T\subseteq \omega^{<\omega}$ be trees. A function $f:S\rightarrow T$ is said to be order preserving if and only if for all
$\sigma,\tau\in S$, if $\sigma\subset\tau$, then $f(\sigma)\subset f(\tau)$. 
\end{enumerate}
\end{definition}
\begin{exercise}

\begin{enumerate}
\item Show that $T\in \mathcal{WF}_\alpha$ iff there is an $<_T$-order preserving function $f:T\to \alpha$.
\item From this  deduce that
 for every $\alpha<\omega_1$, $\mathcal{W}\mathcal{F}_{\alpha+1}\setminus\mathcal{W}\mathcal{F}_\alpha\not=\emptyset$.
 \end{enumerate}
\end{exercise}

\begin{lemma}\label{orderpres}
Let $S,T\subseteq \omega^{<\omega}$ be trees. 
\begin{enumerate}
	\item $\rho(S)\leq\rho(T)$ if and only if there exists an order preserving $f:S\rightarrow T$. Moreover, if $f$ exists and $S\ne \emptyset$ we can assume that $f(\emptyset)=\emptyset$.
	\item If $T$ is well founded, then 
	
\[\rho(S)<\rho(T)\]  
\[\mbox{iff}\]
\[[(S=\emptyset\land T\not=\emptyset)\]
\[  \mbox{ or}\]
\[
\mbox{there is an $n\in \om$, $\la n\ra\in T$ and an order preserving }
 f:S\rightarrow T \mbox{ with }f(\emptyset)=\la n\ra].\]									
\end{enumerate}
\end{lemma}

\pf
	To see $(1)$, suppose first that $f:S\rightarrow T$ is order preserving. Without loss of generality we can assume that $T$ is well-founded. We show by induction on $\alpha$ that if
	$\rho_S(\sigma)=\alpha$, then $\rho_T(f(\sigma))\geq\alpha$. This will suffice as $\rho(S)=\rho_S(\emptyset)\leq
	\rho_T(f(\emptyset))\leq\rho_T(\emptyset)=\rho(T)$.
	
	The case $\alpha=0$ is trivial. Suppose the induction hypothesis holds for all $\beta<\alpha$. 
	Let $\sigma$ be such that $\rho_S(\sigma)=\alpha$. Then by hypothesis
	\begin{eqnarray*}
		\rho_S(\sigma) & = & \sup\{\rho_S(\sigma{\cat}i)+1\ |\ \sigma\cat i\in S\land i\in\omega\}\\
		              & \leq & \sup\{\rho_T(f(\sigma{\cat}i))+1\ |\ \sigma\cat i\in S \land i\in\omega\}\\
					  & \leq & \rho_T(f(\sigma))
					  \end{eqnarray*}
	The last inequality follows from the fact that $\rho_T(f(\sigma{\cat}i))\leq\rho_T(f(\sigma){\cat}i)$ for every $i$ as
	$f(\sigma)\subset f(\sigma{\cat}i)$.
	
	Conversely, suppose that $\rho(S)\leq \rho(T)$. 
	
\medskip
\noindent{\bf Case 1:} $T$ is well-founded.
	Define $f:S\rightarrow T$ by induction on $lh(\sigma)$ for $\sigma\in S$ with the property that for every 
	$\sigma$, $\rho_S(\sigma)\leq\rho_T(f(\sigma))$. If $lh(\sigma)=0$, then $\sigma=\emptyset$, so let $f(\emptyset)=\emptyset$.
	Suppose we have defined $f$ on all sequences of length $n$ in $S$, and fix $\sigma\in S$ with $lh(\sigma)=n$.
	Then $\rho_S(\sigma)=\sup\{\rho_S(\sigma{\cat}i)+1\ |\ i\in\omega\}$ and $\rho_T(f(\sigma))=\sup\{\rho_T(f(\sigma){\cat}j)+1\ |\ j\in\omega\}$.
	By the induction hypothesis, for every $i$ there is a $j_i$ such that $\rho_S(\sigma{\cat}i)\leq\rho_T(f(\sigma){\cat}j_i)$. 
	So let $f(\sigma{\cat}i)=f(\sigma){\cat}j_i$. 
	This function is as desired.
\medskip

\noindent{\bf Case 2:} $T$ is ill-founded. 
Let $b:\omega\to \omega$ be a branch through $T$. Define $f:S\to T$ by $f(\sigma)=\la b(0), b(1), \dots b(lh(\sigma)-1)\ra$.
	
	\medskip
	To prove $(2)$, let $T$ be well founded. Suppose first that $\rho(S)<\rho(T)$ and $S$ is not empty. Then 
\begin{eqnarray*}\rho(T)&=&\rho_T(\emptyset)\\
&=&\sup_{n\in\omega}\{\rho_T({\la}n{\ra})+1\}\\
&=&\sup_{n\in\omega}\{\rho(T(\la n\ra)+1\}.
\end{eqnarray*}
	So there is an $n\in\omega$ such that $\rho(S)\leq\rho(T({{\la}n{\ra}}))$. By $(1)$, construct an order preserving $f_0: S\rightarrow T({{\la}n{\ra}})$ with $f_0(\emptyset)=\emptyset$. The $f_0$ induces an order preserving $f:S\to T$ with $f(\emptyset)=\la n \ra$.
	
	Conversely, suppose $f: S\rightarrow T({{\la}n{\ra}})$ is order preserving. Then $\rho(S)\leq \rho(T({{\la}n{\ra}}))<\rho(T)$ because $T$ is well founded.
\qed
\begin{lemma}
\(\mathcal{W}\mathcal{F}_\alpha\) is Borel.
\end{lemma}
\pf
We proceed by induction on $\alpha$.
Since $\emptyset$ is Borel, the case $\alpha=1$ is clear. For any ordinal $\alpha<\omega_1$, we have
\[\mathcal{W}\mathcal{F}_{\alpha+1}=\bigcap_{n\in\omega}\{ T\ |\ {\la}n\ra\not\in T \lor T({{\la}n\ra})\in\mathcal{W}\mathcal{F}_{\alpha}\}\]
which is Borel.

Suppose $\alpha$ is a limit ordinal. Then $\mathcal{W}\mathcal{F}_{\alpha}=\bigcup_{\beta<\alpha}\mathcal{W}\mathcal{F}_{\beta}$ which is also Borel.
\qed

\begin{theorem}
Every $\bpi_1^1$ set is a union of $\omega_1$ many Borel sets.
\end{theorem}
\pf
	Let $A\subseteq X$ be a $\bpi_1^1$ set. Let $f:X\rightarrow \trees$ be a Borel reduction of $A$ to the well founded trees. 
	Then $A=\bigcup_{\alpha\in\omega_1} f^{-1}(\mathcal{W}\mathcal{F}_{\alpha})$.\qed

\begin{exercise}
Show that if $A\subseteq \omom$ is a $\bpi^1_1$ set and $|A|>\omega_1$ then $|A|=2^{\aleph_0}$.
\end{exercise}	

\begin{exercise}
Show that every analytic subset of $\omom$ is a union of $\omega_1$ Borel sets.

\noindent({\bf Hint:} Suppose that $A=p[T]$. Let $B_\alpha=\{x:(\rho(T^x)\not<\alpha) \land (\forall s\in \om^{<\om})(\rho(T^x(s))\ne \alpha)\}$. Show that $B_\alpha$ is Borel and $A=\bigcup B_\alpha$.)
\end{exercise}

\begin{lemma} \label{relationlemma}
\begin{enumerate}
\item $\{(S,T)\subseteq \omega^{<\omega}\times\omega^{<\omega}\ |\ \rho(S)\leq\rho(T)\}$ is a $\bsigma_1^1$ subset of the Cartesian product of the space of 
trees with itself.
\item There is a $\bsigma_1^1$ set $R\subseteq \trees\times\trees$ such that if $T$ is well founded, then 
\[\{S\ |\ (S,T)\in R\} = \{S\ |\ \rho(S)<\rho(T)\}\]
\end{enumerate}
\end{lemma}

\pf
	For $(1)$, we have that $\rho(S)\leq\rho(T)$ if and only if there is a $f:\omega^{<\omega}\rightarrow\omega^{<\omega}$ such that $f$ is order
	preserving from $S$ to $T$. In other words, 
	\[\rho(S)\leq\rho(T) \]
	iff
	\[(\exists F\in 2^{{\omega^{<\omega}}\times\omega^{<\omega}})((F\mbox{ is a function with domain $S$) and
	($F$ is order preserving) and}\]
\[ \forall\sigma(\sigma\in S\Rightarrow F(\sigma)\in T))\]
	That $F$ is a function is a closed condition, that it is order preserving is a closed condition, and for each $\sigma$,  $\sigma\in S\Rightarrow F(\sigma)\in T$
	is  a clopen condition. 
	
	For $(2)$, let $R=\{(S,T)\ |\ (\exists f\in 2^{\omega^{<\omega}\times\omega^{<\omega}})(\exists n\in\omega)( f\mbox{ is order
	preserving from $S$ to $T({{\la}n}\ra)$})\}$. By the second part of Lemma \ref{orderpres}, $R$ is as desired.
\qed
\begin{theorem}[Boundedness Theorem]\label{boundedness theorem}
	Suppose $A\subseteq \mathcal{W}\mathcal{F}$ is $\bsigma^1_1$. Then there is an $\alpha<\omega_1$ such that $A\subseteq \mathcal{W}\mathcal{F}_\alpha$.
	
\end{theorem}
\pf
	Suppose not. Then $S$ is well founded if and only if $(\exists T\in A)(\rho(S)\leq\rho(T))$. By the previous lemma, the set of 
	well founded trees is a $\bsigma_1^1$ subset of $\trees$. Since the set of well founded trees is $\bpi_1^1$, this shows they form a Borel set, a contradiction.\qed

\begin{exercise}
Show that for $A\subseteq {\mathcal LO}$, if $A\subseteq {\mathcal WO}$ is analytic then  there is an $\alpha\in\omega_1$ such that for all $L\in A$ the order-type of $L$ is less than $\alpha$.
\end{exercise}
%





\begin{corollary}\label{boundedness for sets of reals}
Suppose that $R\subseteq \omom\times \omom$  is $\bsigma^1_1$ and $R$ is well-founded. Then $R$ has height less than $\omega_1$.
\end{corollary}
\pf
We start with an abstract claim:
\medskip

\noindent{\bf Claim:} Suppose that $R$ is an arbitrary well-founded relation on a set $X$ of height at least $\omega_1$. Then for all countable linear orderings $I=(B, <_I)$:
\[<_I \mbox{ is a well-ordering}\]
\begin{center} iff
\end{center}
\smallskip
there is a countable set $S\subseteq X$ and a surjection $f:S\to B$ such that:
\begin{quotation}
\noindent for all $s\in S$ and for all $b<_If(s)$ there is a $t\in S$ with $tRs$ and $b\le_I f(t).$
\end{quotation}

To see the claim: Suppose that $<_I$ is a well-ordering of length $\alpha<\omega_1$. Let $x_0\in X$ be such that $ht_R(x_0)=\alpha$. Build a sequence of countable sets $\la S_n:n\in\om\ra$ by induction so that:
\begin{enumerate}
\item $S_0=\{x_0\}$,
\item if $y\in S_n$ there are $\{z_i:i\in\om\}\subseteq S_{n+1}$ such that for all $i, z_i R y$ and $\sup\{ht_R(z_i)+1:i\in\om\}= ht_R(y)$.
\end{enumerate}
Let $S=\bigcup_n S_n$. Then $\{ht_{R\rest S}(s):s\in S\}=\alpha$. Let $f:S\to B$ be defined by setting 
\[f(z)=m \mbox{ iff } ht_R(z)=ht_I(m).\]

For the other direction: towards a contradiction, suppose that $\la b_n:n\in\om\ra$ is a $<_I$-decreasing sequence and $f$ is a function with the stated property. Inductively choose a sequence $\la s_n:n\in \om\ra$ such that
\[s_{n+1}Rs_n \mbox{ and }b_{n+1}\le f(s_{n+1}).\] 
The existence of such a sequence contradicts the well-foundedness of $R$. 

Using the claim we can give a $\bsigma^1_1$ definition of ${\mathcal WO}$. If $I\in {\mathcal LO}$, we have:

\[ I\in {\mathcal WO} \mbox{ iff } (\exists \la z^i_n:i,n\in\om\ra\in(\omom)^{\om\times \om})\]
\[(\forall m)(\forall n)(m<_In\implies(\forall i)(\exists k)(\exists j)(z^j_kRz^i_n\land m\le_Ik)).\]
\qed

\subsection{ Norms}
In this section we formalize the properties we have shown about the well-founded sets. 

\begin{definition}
\
\begin{itemize}
\item A norm on a set $A\subseteq X$ is a function $\phi: X\rightarrow{\mathcal{ON}}\cup\{\infty\}$, the class of ordinals such that $x\in A$ iff $\phi(x)\ne \infty$. The set $\{(x, y):\phi(x)\le \phi(y)\}$ is  called a \emph{pre-well-ordering} of $A$. We treat ``$\infty$" as a symbol that is bigger than all of the ordinals. Since $A$ is a set, we could also take $\infty$ to be any particular ordinal larger than $\sup  \phi[A]$.

It is an easy set-theoretic exercise to show that for every norm $\phi$ there is a norm $\psi$ such that $\phi$ and $\psi$ have the same pre-well-orderings and there is an ordinal $\beta$ such that $\psi:A\to \beta$ is surjective.

\item Let $A\in\bpi_1^1(X)$, and $\phi:A\rightarrow{\mathcal{ON}}$ be a norm. Then $\phi$ is a $\bpi_1^1$-norm if and only if 
 there are relations $\leq^{\bpi^1_1}\in\bpi_1^1(X\times X)$ and $\leq^{\bsigma^1_1}\in\bsigma_1^1(X\times X)$ such that for all $y\in A$
 \begin{eqnarray*}
 (x\in A)\land (\phi(x)\leq\phi(y)) & \Leftrightarrow & x\leq^{\bpi^1_1}y\\
								& \Leftrightarrow & x\leq^{\bsigma^1_1}y
\end{eqnarray*}
\end{itemize}
\end{definition}

Here are an equivalent forms of this definition: 
\begin{prop}Let $X$ be Polish, $A$ be a $\bpi^1_1$-set and $\phi:X\to \mathcal{ON}\cup \{\infty\}$ be such that $x\in A$ iff $\phi(x)<\infty$. Then $\phi$ is a $\bpi^1_1$-norm iff 
the relations $\le^*_\phi$ and $<^*_\phi$ are both $\bpi^1_1$ where:
\begin{enumerate}
\item $x\le^*_\phi y$ iff $x\in A \land (\phi(x)\le \phi(y))$,
\item $x<^*_\phi y$ iff $x\in A \land \phi(x)<\phi(y)$.
\end{enumerate}
\end{prop}

\begin{theorem}
The function $\rho:\mathcal{W}\mathcal{F}\rightarrow \omega_1$ is a $\bpi_1^1$-norm.
\end{theorem}
\pf
	The relation $\leq^{\bsigma^1_1}$ is the one in part $(1)$ of Lemma \ref{relationlemma}. 
	
	Let $R$ be the relation given in part $(2)$ of this same Lemma. Then 
	\[x\leq^{\bpi^1_1}y\Leftrightarrow (y,x)\not\in R\land x\in \mathcal{W}\mathcal{F}\]
	so that $\leq^{\bpi^1_1}$ is the intersection of two $\bpi_1^1$ sets.

\qed
\begin{corollary}
	Let $A\in\bpi_1^1(X)$, and $f:X\rightarrow{\trees}$ be a Borel reduction of $A$ to the well founded trees.
	Then $\rho\circ f: A\rightarrow \omega_1$ is a $\bpi_1^1$-norm.

\end{corollary}
\pf
	Let $\leq^{\bpi^1_1}_{{\trees}}$ and $\leq^{\bsigma^1_1}_{{\trees}}$ be witnesses to $\rho$ being a $\bpi_1^1$-norm on well founded trees.
	Define $\leq^{\bpi^1_1}=\{(x,y)\ |\ f(x)\leq^{\bpi^1_1}_{{\trees}}f(y)\}$ and
	$\leq^{\bsigma^1_1}=\{(x,y)\ |\ f(x)\leq^{\bpi^1_1}_{{\trees}}f(y)\}$ which are the witnesses for $\rho\circ f: A\rightarrow \omega_1$ 
	being a $\bpi_1^1$-norm.
\qed

\begin{prop}
Suppose that $A\subseteq X$ is a $\bpi^1_1$ set, with $X$ Polish and that $\phi:A\to \omega_1\cup\{\infty\}$ is a $\bpi^1_1$-norm. Then $A$ is Borel iff there is a countable ordinal $\alpha$ such that $\phi[A]\subseteq \alpha$.
\end{prop}
\pf
Suppose that $A$ is Borel. Since $A$ is $\bsigma^1_1$ the following relation is $\bsigma^1_1$:
\[xRy \mbox{ iff } x\in A\land y\in A \land x\le^{\bsigma^1_1} y.\]
By Corollary \ref{boundedness for sets of reals} we know that the height of $R$ is countable. Hence there is an $\alpha$ such that the range of $\phi$ is a subset of $\alpha$. 

For the other direction: suppose that there is a countable ordinal $\alpha$ such that the range of $\phi$ is a subset of $\alpha$. Then we can find a sequence $\la z_i:i\in\om\ra$ of elements of $A$ such that $\{\phi(z_i):i\in\om\}$ is cofinal in $\sup \phi[A]$. Then $A$ is defined by:
\[x\in A\mbox{ iff }(\exists i)(x\le^{\bsigma^1_1} z_i).\]
Hence $A$ is $\bsigma^1_1$. 
\qed

\subsection{Reduction and Separation}

\begin{definition}
	A class $\Gamma$ of sets has the \emph{reduction property} if and only if for any Polish space $X$ and for all $A,B\in\Gamma(X)$, there are $A_0\subseteq A$, $B_0\subseteq B$ with:
	\begin{itemize}
	\item $A_0,B_0\in\Gamma(X)$
	\item $A_0\cap B_0=\emptyset$
	\item $A_0\cup B_0=A\cup B$
	\end{itemize}
\end{definition}

\begin{theorem}
	Let $X$ be a perfect Polish space. Then $\bpi_1^1(X)$ has the reduction property.

\end{theorem}
\pf
	Let $f,g$ be reductions of $A$, $B$ to $\mathcal{W}\mathcal{F}$. Let 
	$A_0=\{x\in A\ |\ \rho(f(x))\leq \rho(g(x))\}$ and $B_0=\{x\in B\ |\ \rho(f(x))> \rho(g(x))\}$.
	It is clear that $B_0$ is $\bpi_1^1$ as it is the complement of a $\bsigma_1^1$ set by Lemma \ref{relationlemma}. 
	For every $x\in A$, we have that $f(x)$ is well founded, so that $A=\{x\in A\ |\ \rho(f(x))\leq^{\bpi_1^1} \rho(g(x))\}$ is also $\bpi_1^1$.
	
	Finally, notice that if $x\in A\setminus A_0$, then in fact $x\in B_0$ (the reverse is also true). Hence $A_0\cup B_0= A\cup B$.

\qed

Notice that if $C$ and $D$ are analytic, we can apply the reduction property to their complements to obtain a 
$\bdelta_1^1$ (hence Borel) set $S$ separating $C$ and $D$.
In general, the reduction property for $\Gamma$ implies the separation for its dual class.

\begin{proposition}
	$\bpi_1^1$ does not have the separation property.
	
\end{proposition}

\pf
	Let $U\in\bpi_1^1({\omega^\omega\times\omega^\omega})$ be universal for $\bpi_1^1(\omega^\omega)$. Let $h:\omega^\omega\rightarrow\omega^\omega\times\omega^\omega$ be a homeomorphism. Denote  
	$h(x)=(x_{0},x_{1})$.
	
	Let $P=\{x\ |\ (x_{even},x)\in U\}$ and $Q=\{x\ |\ (x_{odd},x)\in U\}$. Using $\bpi_1^1$ reduction, produce $P_0$ and $Q_0$ that reduce $P$ and $Q$. 
	
	Let us assume, for a contradiction, that there is a Borel set $B$ with $B\cap Q_0=\emptyset$ and $P_0\subseteq B$.
	By universality of $U$, there are $a,b\in\omega^\omega$ such that $B=U_a$ and $\omega^\omega\setminus B=U_b$.
	Let $x=h^{-1}(b,a)$.  Note that if $x\in Q\cap B$ then $x\notin Q_0$ so $x\in Q\cap P$.
	\begin{eqnarray*}
		x\in B & \Leftrightarrow & (a,x)\in U\\
				& \Leftrightarrow & x\in Q\cap P\\
				& \Leftrightarrow & (b,x)\in U\\
				& \Leftrightarrow & x\in U_b\\
				& \Leftrightarrow & x\in\omega^\omega\setminus B
	\end{eqnarray*}
	This contradiction proves the proposition.\qed

\section{Uniformization}

\begin{definition}
	Let $X$, $Y$ be Polish, $A\subseteq X\times Y$. We say that a set $B$ uniformizes $A$ if and only if:
	\begin{enumerate}
	\item $\Pi_X(B)=\Pi_X(A)$
	\item for every $x\in\Pi_X(A)$ there is a unique $y\in Y$ such that $(x,y)\in B$
	\end{enumerate}
\end{definition}

Clearly $B$ can be viewed as a choice function defined on $\Pi_X(A)$ and thus the axiom of choice guarantees the existence of a uniformization for any set. However, we are only interested in ``nice" uniformizations. 

We make this into a definition:

\begin{definition}
Let $X, Y$ be Polish and $A\subseteq X\times Y$. A function $f:A\to Y$ \emph{uniformizes} $A$ iff $B=$graph$(f)$ uniformizes $A$.
\end{definition}

\subsection{The bad news}
\begin{example}
There is a closed set $C\subseteq \omega^\omega\times\omega^\omega$ that cannot be uniformized by any analytic set.
\end{example}
\pf
	Let $A_0, A_1\in\bpi_1^1(\omega^\omega)$ be disjoint and inseparable by any Borel set. Then there are closed sets $C_0, C_1$ such that $\omom\setminus A_i=\{x:(\exists y)(x,y)\in C_i\}$. Without loss of generality we can assume that $C_i\subseteq \omom\times [\la i\ra]$, where $[\la i\ra]$ consists of those elements of $\omom$ such that $y(0)=i$. 
(Exercise: why can we assume this?)
 	
Let $C=C_0\cup C_1$. Then $C$ is  a closed subset of $\omom$, and since $A_0\cap A_1=\emptyset$, $\Pi_X(C)=\omom$.  Suppose that $C$ can be uniformized by an analytic set $A$. Since the projection of $C$ is $\omom$, the analytic set $A$ is the graph of an analytic function $f:\omom\to \omom\times \omom$. Hence, by Corollary \ref{analytic functions are Borel}, $f$ is Borel measurable. 

Let $B_i=f^{-1}(\omom\times [\la i\ra])$ for $i=1, 2$. Then $B_i\subseteq \omom\setminus A_i$, $B_0\cap B_1=\emptyset$,  and $B_0\cup B_1=\omom$. Hence $B_i\supseteq A_{1-i}$ and thus $B_i$ separates $A_{1-i}$ from $A_i$,  a contradiction.	
\qed
\subsection{Another digression on trees}

Suppose that $T\subseteq A^{<\om}$ is  an ill-founded tree and $<$ is a well-ordering of $A$. Then the \emph{left-most branch} through $T$ is the function $b:\om\to A$ with the property that for all $n\in\om, b(n)$ is the $<$-least element of $A$ among $\{c(n):c\rest n=b\rest n\land c\in [T]\}$.

We note that there is always a unique left most branch and the map ``$T\mapsto$ the leftmost branch of T", is a continuous map from the closed set $\{T\subseteq A^{<\om}:T$ is well-pruned$\}\subseteq \trees$ to $A^\om$.

\subsection{Jankov-von Neumann uniformization}
The easiest (and most commonly used) uniformization theorem is often called the ``Jankov-von Neumann" uniformization theorem. Let $\sigma(\bsigma^1_1)$ be the $\sigma$-algebra generated by the $\bsigma^1_1$-sets.  Since every $\bsigma^1_1$-set is universally measurable (Theorem \ref {analytic sets are measurable}) every  function that is measurable with respect to the $\sigma$-algebra  $\sigma(\bsigma^1_1)$ is universally measurable, and hence the next theorem shows that every analytic set can be uniformized by a function that is universally measurable.
\begin{theorem}
Suppose that $X,Y$ are uncountable Polish spaces and $A\subseteq X\times Y$ is analytic. Then  there is a function $f$ that uniformizes $A$ that is measurable with respect to $\sigma(\bsigma^1_1)$.
\end{theorem}
 \pf Since this theorem only refers to the Borel structure of $X$ and $Y$ we can Exercise \ref{borel bijections two} to see that  it holds for all standard Borel spaces. Moreover it suffices to prove the result in the case that $X=Y=\omom$.
 
 Let $T\subseteq (\om\times \om\times\om)^{<\om}$ be a well-pruned tree such that  $A=\{(x,y):T^{x,y}$ is ill-founded$\}$.   Let $A'=[T]$. Then $A'$ is a closed subset of $\omom\times \omom\times \omom$. Call the three copies of $\omom$ by $X, Y$ and $Z$. 
 
 Suppose that we can uniformize $A'\subseteq X \times (Y\times Z)$ by a function $f$ that is measurable with respect to $\sigma(\bsigma^1_1)$. Let $\Pi_Y:Y\times Z\to Y$ be the projection map. Then $\Pi_Y$ is continuous and hence $\Pi_Y\circ f:X\to Y$ is measurable with respect to $\sigma(\bsigma^1_1)$ and hence satisfies the theorem. 
 
 Thus we are reduced to finding the function $f$. Note that $Y\times Z\cong \omom$ so we can assume that we are trying to uniformize a closed subset $A$ of $\omom\times \omom$.
 
 Now let $T\subset (\omom\times \omom)^{<\om}$ be a well-pruned tree  such that $A=[T]$. For each $x\in \Pi_X(A)$, let $f(x)$ be the left-most branch of $[T^x]$. We claim that $f:A\to \omom$ is $\sigma(\bsigma^1_1)$-measurable; i.e. for all $s\in \om^{<\om}$, the set $f^{-1}([s])\in\sigma(\bsigma^1_1)$. We show this by induction on the length of $s$.
 
 Since $A$ is $\bsigma^1_1$, $\Pi_X(A)$ is $\bsigma^1_1$ and hence  this holds if $s=\emptyset$. Suppose that $f^{-1}([s])\in \sigma(\bsigma^1_1)$, we show it for $f^{-1}([s\cat n])$.  But  $x\in f^{-1}([s\cat n])$ iff $x\in f^{-1}([s])$ and 
 \begin{equation}
(\exists y\in [s\cat n])((x,y)\in T)\land (\forall m<n)( T^x(s\cat m)\mbox{ is well-founded}).
 \end{equation}
 
 This set is the intersection of $f^{-s}([s])$ with a $\bsigma^1_1$ set and an intersection of countably many $\bpi^1_1$ sets. Hence it is in $\sigma(\bsigma^1_1)$.\qed
  
 \begin{exercise}
 Show that there is a function $f$ that uniformizes $P$ that is measurable with respect to the $\sigma$-algebra of sets with the property of Baire.
 \end{exercise}

\end{document}